\documentclass[pdflatex,sn-mathphys-num]{sn-jnl}
\usepackage[T1]{fontenc}
\usepackage{lmodern}
\usepackage{amsmath,amssymb,amsthm,mathtools}
\usepackage{booktabs,longtable,array,enumitem,needspace,xurl,xcolor}
\usepackage[title]{appendix}
\hypersetup{hidelinks,
  pdftitle={Positivity and uniform log-concavity for rank-two Chern coefficients},
  pdfauthor={Hang-Cheng Dong and Changle Zhao},
  pdfsubject={Chern classes, Schur coefficients and log-concavity}}
\setlist[enumerate]{label=\textup{(\roman*)},leftmargin=2em,itemsep=3pt}
\setlist[itemize]{leftmargin=1.7em,itemsep=3pt}
\theoremstyle{thmstyleone}
\newtheorem{theorem}{Theorem}[section]
\newtheorem{proposition}[theorem]{Proposition}
\newtheorem{lemma}[theorem]{Lemma}
\newtheorem{corollary}[theorem]{Corollary}
\theoremstyle{thmstylethree}

\newtheorem{conjecture}[theorem]{Conjecture}
\theoremstyle{thmstyletwo}
\newtheorem{remark}[theorem]{Remark}

\numberwithin{equation}{section}
\newcommand{\Ccal}{\mathcal C}
\newcommand{\F}{\mathcal F}
\newcommand{\fall}[2]{(#1)_{\underline{#2}}}

\DeclareMathOperator{\supp}{supp}
\DeclareMathOperator{\Sym}{Sym}

\newcolumntype{P}[1]{>{\raggedright\arraybackslash}p{#1}}
\begin{document}

\title[Rank-two Chern coefficients]{Positivity and uniform log-concavity
for rank-two Chern coefficients}

\author[1]{\fnm{Hang-Cheng} \sur{Dong}}\email{hunsen\_d@163.com}
\author[2]{\fnm{Changle} \sur{Zhao}}\email{646627477@qq.com}
\affil[1]{\orgdiv{School of Instrumentation Science and Engineering},
  \orgname{Harbin Institute of Technology},
  \orgaddress{\city{Harbin}, \country{China}}}
\affil[2]{\orgname{Hangzhou Mailiang Technology Co., Ltd.},
  \orgaddress{\city{Hangzhou}, \country{China}}}

\abstract{We study the binomial-basis coefficients in the Schur expansion
of the Chern classes of symmetric powers of a rank-two bundle.
B\'erczi and Feh\'er conjectured that these coefficients are nonnegative
and log-concave in the Newton index.  An inverse-series formula yields
an integral gamma-positive refinement with exact support and a positive
expression for the derivatives of the shifted Chern quotient.
We prove that generalized Bernoulli polynomials of real order at least
their degree have simple real zeros.  Applied to the shifted quotient,
this gives strict ultra-log-concavity in the Schur index at every positive
shift.  For the original Newton coefficients, we establish strict
log-concavity with an explicit multiplicative margin in two uniform
regions.  A Hahn projection gives a square-root band about the balanced
Schur index; a complex saddle estimate extends this to a band of linear
width.  Uniform estimates for the cumulants of a normalized Schur
functional give a terminal Newton region of linear width.
Polynomial coefficient identities also establish every row with Schur
index at most six or imbalance at most nine.  The argument separates
log-concavity in the Schur, derivative and Newton indices and identifies
the factorial factors in the transfer between the latter two.
The unrestricted Newton log-concavity conjecture remains open.}

\maketitle

\section{Introduction}
\label{sec:introduction}

Let $E$ be a vector bundle of rank two with formal Chern roots $x,y$.
The weights of $\Sym^d E$ are $ix+(d-i)y$, $0\le i\le d$.  Its
$k$th Chern class is therefore represented by
\begin{equation}\label{intro:chern}
 C_{k,d}(x,y)=[z^k]\prod_{i=0}^d
   \bigl(1+z(ix+(d-i)y)\bigr)
 =\sum_{j=0}^{\lfloor k/2\rfloor}A_{k,j}(d)s_{(k-j,j)}(x,y),
\end{equation}
where $s_{(k-j,j)}$ denotes a two-variable Schur polynomial.
For fixed $k,j$, the coefficient $A_{k,j}(d)$ extends to a polynomial
in $d$ of degree at most $2k$.  We write its two Newton expansions as
\begin{equation}\label{intro:newton}
 \begin{gathered}
 A_{k,j}(d)=\sum_N B_{k,j,N}\binom dN
          =\sum_N g_{k,j,N}\binom{d+1}{N},\\
 B_{k,j,N}=g_{k,j,N}+g_{k,j,N+1}.
 \end{gathered}
\end{equation}
Sequences are extended by zero outside their support.  We call a
nonnegative sequence log-concave if its support is an interval and
$c_i^2\ge c_{i-1}c_{i+1}$ at every index.

B\'erczi and Feh\'er \cite[Conjectures~3.2 and~3.3]{BF} conjectured
that each row $(B_{k,j,N})_N$ is nonnegative and log-concave.  They
proved log-concavity for $j=0,1$ \cite[Theorems~3.7 and~5.4]{BF}.
A subsequent argument for nonnegativity at every Schur index was
posted by Patel \cite{Patel}.  The expansion in the shifted basis
$\binom{d+1}{N}$ reveals a stronger positivity property: each coefficient
has a nonnegative integral expansion in $(xy)^p(x+y)^{k-2p}$.
We determine its support and use the same inverse series to study
three associated log-concavity questions.

Set $\sigma=x+y$ and $\rho=xy$, and define
\begin{equation}\label{intro:inverse}
 \phi(t)=\frac{e^{-yt}-e^{-xt}}{x-y},\qquad T=\phi^{-1},\qquad
 \Ccal=\log T'-\sigma T,\qquad H=\log(T/z).
\end{equation}
The quotient defining $\phi$ has a removable specialization at $x=y$.
Since $\phi(t)=t+O(t^2)$, its compositional inverse exists.

\begin{theorem}[Integral gamma refinement]\label{intro:gamma}
For $k\ge1$, there are nonnegative integers $\Gamma_{k,p,N}$ such that
\begin{equation}\label{intro:gamma-expansion}
 C_{k,d}(x,y)=\sum_{p,N}\Gamma_{k,p,N}\binom{d+1}{N}
                  \rho^p\sigma^{k-2p}.
\end{equation}
Their positive support is exactly
\[
 0\le p\le\lfloor k/2\rfloor,\qquad
 k\le N\le2k-p,\qquad p>0\ \text{or}\ N>k.
\]
Consequently,
\[
 \supp g_{k,j,\bullet}=
 \begin{cases}[k+1,2k],&j=0,\\{}[k,2k],&j\ge1,\end{cases}
 \qquad
 \supp B_{k,j,\bullet}=
 \begin{cases}[k,2k],&j=0,\\{}[k-1,2k],&j\ge1.\end{cases}
\]
\end{theorem}

The proof follows from the differential equations
\[
 V'=\sigma V^2+\rho zV^3,\qquad \Ccal'=\rho zV^2,
 \qquad V=T'.
\]
Lagrange inversion also identifies the quotient by the forced factorial:
\begin{equation}\label{intro:quotient}
 A_{k,j}(d)=\fall{d+1}{k}P_{k,j}(d),\qquad
 F_{k,j}(s)=P_{k,j}(s+k-1).
\end{equation}
Here $\fall{u}{k}=u(u-1)\cdots(u-k+1)$.  The shifted quotients satisfy
\begin{equation}\label{intro:shifted}
 \F_k(s;x,y):=\sum_jF_{k,j}(s)s_{(k-j,j)}(x,y)
       =[z^k]\exp(\Ccal+sH).
\end{equation}
Their derivatives $a_{k,j,r}:=F_{k,j}^{(r)}(0)$ are positive for
$1\le r\le k$, while $a_{k,j,0}>0$ precisely when $j\ge1$.
Section~\ref{src:section} derives these assertions from
\eqref{intro:chern}.

Real-rootedness in the Chern variables provides a second consequence
of the quotient formula.  With the normalization
\[
 B_n^{(\alpha)}(w)=n![t^n]
       \left(\frac{t}{e^t-1}\right)^\alpha e^{wt},
\]
the required Bernoulli root theorem holds for real order.

\begin{theorem}\label{intro:bernoulli}
For every integer $n\ge1$ and real $\alpha\ge n$,
$B_n^{(\alpha)}$ has $n$ simple real roots in $(0,\alpha)$.
If $\alpha>n+1$, all roots lie in $(1,\alpha-1)$, and
$B_n^{(n+1)}(w)=\prod_{i=1}^n(w-i)$.
For every real $s>0$, the polynomial $\F_k(s;x,1)$ has $k$ simple
negative roots.  Moreover, for $1\le j<\lfloor k/2\rfloor$,
\begin{equation}\label{intro:schur-ulc}
 F_{k,j}(s)^2>
 \frac{(j+1)(k+2-j)}{j(k+1-j)}F_{k,j-1}(s)F_{k,j+1}(s).
\end{equation}
\end{theorem}

Thus $F_{k,j}(s)/\binom{k+1}{j}$ is strictly log-concave in $j$ at
each positive shift.  The Bernoulli theorem extends the integer-order
root mechanism discussed by \v Stampach \cite{Stampach}; its proof
uses a cut-contour sign calculation followed by Rolle's theorem.
The inequality in $j$ is distinct from log-concavity of the derivative
row in $r$ and of the Newton row in $N$.  For the latter problem, the
derivative row provides a sufficient condition.  The transfer theorem
in Section~\ref{sec:bridge} retains the factorial contribution and
yields the factor
\begin{equation}\label{intro:factor}
 \beta_{k,N}=1+\frac{k}{(N-k+1)(N+2)},\qquad k\le N\le2k-1.
\end{equation}

Write $k=2n+\epsilon$, $j=n-m$, where $\epsilon\in\{0,1\}$ and
$0\le m\le n$, and put $\omega=2m+\epsilon+1$.

\begin{theorem}[Uniform Newton inequalities]\label{intro:regions}
For every admissible Schur index, the following assertions hold.
\begin{enumerate}
\item If $m\ge2$ and $n\ge5(m+2)^2$, then
\begin{equation}\label{intro:strict}
 B_{k,j,N}^2>\beta_{k,N}B_{k,j,N-1}B_{k,j,N+1}
 \qquad(k\le N\le2k-1).
\end{equation}
The same conclusion holds if $n\ge10^{40}$ and
$0<\omega\le10^{-12}n$.
\item Let $Q\ge1$.  If $k\ge2^{38}(Q+1)$, then
\eqref{intro:strict} holds for every $j$ at all centers
$2k-Q\le N\le2k-1$.
\end{enumerate}
\end{theorem}

The first part concerns full Newton rows near the balanced Schur
index; the second concerns the terminal part of every row.
The constants are sufficient bounds and are not optimized.
They imply regions of linear width in both directions.

\begin{corollary}\label{intro:linear-regions}
Set $c=10^{-13}$ and $K_0=3\cdot10^{40}$.  If $k\ge K_0$, then
\eqref{intro:strict} holds whenever
\[
 k-2j\le ck\qquad\text{or}\qquad 2k-N\le ck,
 \qquad k\le N\le2k-1.
\]
\end{corollary}
\begin{proof}
In the first case, $n=\lfloor k/2\rfloor\ge10^{40}$ and
$\omega\le10^{-13}(2n+\epsilon)+1\le10^{-12}n$.
In the second, take $Q=2k-N\ge1$.  Then
$2^{38}(Q+1)\le2^{38}(10^{-13}k+1)<k$.
Theorem~\ref{intro:regions} applies in each case.
\end{proof}

The remaining low-index families are obtained from polynomial
coefficient identities and finite exact inequalities.

\begin{theorem}[Boundary families]\label{intro:families}
For all admissible $k,j$, the row $(B_{k,j,N})_N$ is log-concave if
\[
 j\le6\qquad\text{or}\qquad k-2j\le9.
\]
\end{theorem}

The cases $j=0,1$ are due to B\'erczi and Feh\'er.
Section~\ref{sec:computational} and Appendix~\ref{supp:fixed-section}
treat the other cases.  This theorem asserts ordinary log-concavity.
The stronger factor in \eqref{intro:strict} applies to the subfamilies
for which the derivative inequalities are established.

We briefly describe the two estimates underlying
Theorem~\ref{intro:regions}.  Near the balanced index, Schur extraction
becomes a finite Hahn transform of a hyperbolic-sine generating function.
A zero-disk bound for $P_{k,j}$ places the zeros of $F_{k,j}$ in a
sector of half-angle $\pi/4$ about the negative real axis.  This implies
derivative log-concavity.  Direct coefficient bounds give the quadratic
condition on $n$; a two-variable saddle estimate, uniform in a complex
frequency parameter, gives the linear range.  Tang
\cite[Appendix~B.2]{Tang} discusses a related Newton root-region
criterion.  The sector argument used here is proved in
Section~\ref{sec:near}.

For the terminal Newton indices, the normalized Schur functional has
moments rational in two independent parameters.  Connected expansions
reduce the degree of each logarithmic coefficient in the defect.
The relation between approximate factorization and small cumulants
appears in Do\l\k{e}ga and F\'eray \cite{DF}.  We obtain uniform bounds
for the Schur functional in this problem; a substitution into the
positive generating series then gives a nonvanishing disk of radius
inversely proportional to the defect.  Cauchy estimates in that disk
control the discrete curvature of the derivative sequence.

The passage from derivatives to forward differences uses the
convolution-array principle of F\"oldes and Major \cite{FM}.
For fixed Schur index, the polynomial kernels satisfy the centered-prefix
criterion for LC-positivity of Wang and Yeh \cite{WY}.
Gamma positivity and asymptotic log-concavity also arise in related
enumerative-geometric settings; see Aluffi, Chen and Marcolli \cite{ACM}.

Section~\ref{src:section} establishes the positive inverse-series
formulas, and Section~\ref{bern:section} proves the Bernoulli root theorem
and its Schur consequence.  Section~\ref{sec:bridge} gives the Newton
transfer.  The uniform estimates occupy Sections~\ref{sec:near}
and~\ref{sec:tail}, with their longer analytic arguments in
Appendices~\ref{near:linear-proof} and~\ref{tail:appendix}.
Section~\ref{sec:computational} treats the boundary families, and
Section~\ref{sec:discussion} states the remaining conjecture.

\section{An inverse source and integral gamma positivity}\label{src:section}

For $k\ge1$, let
\begin{equation}\label{src:Chern}
 \begin{aligned}
 C_{k,d}(x,y)&=[u^k]\prod_{i=0}^{d}
       \bigl(1+u(ix+(d-i)y)\bigr)\\
 &=\sum_{j=0}^{\lfloor k/2\rfloor}
       A_{k,j}(d)s_{(k-j,j)}(x,y),
       \qquad d\in\mathbb Z_{\ge0}.
 \end{aligned}
\end{equation}
A formal inverse gives the polynomial continuation in $d$ and the two
Newton expansions
\begin{equation}\label{src:Newton-bases}
 \begin{aligned}
 A_{k,j}(d)&=\sum_{N\ge0}g_{k,j,N}\binom{d+1}{N}
          =\sum_{N\ge0}B_{k,j,N}\binom dN,\\
 W_{k,N}(x,y)&=\sum_j g_{k,j,N}s_{(k-j,j)}(x,y).
 \end{aligned}
\end{equation}
The coefficients $B_{k,j,N}$ occur in the positivity and log-concavity
conjectures of \cite{BF}.  Patel \cite{Patel} has posted an argument
for nonnegativity at every Schur index.  We refine this positivity in the basis
$\binom{d+1}{N}$ and determine the joint support of the refinement.

Throughout this section,
\[
 \sigma=x+y,\qquad \rho=xy,\qquad \delta=x-y.
\]
Define formal series over $\mathbb Q[x,y]$ by
\begin{equation}\label{src:inverse}
 \phi(t)=\frac{e^{-yt}-e^{-xt}}{x-y},\qquad
 \phi(T(z))=z,\qquad V=T',\qquad
 \Ccal=\log V-\sigma T.
\end{equation}
The coefficients of $\phi$ are polynomials in $x,y$; in particular,
the specialization $x=y$ is well defined.  Since $\phi(t)=t+O(t^2)$,
its inverse satisfies
$T(0)=0$, $V(0)=1$, and $\Ccal(0)=0$.

\subsection{Identification with the original product}

\begin{lemma}\label{src:product-source}
For integers $n\ge k\ge0$,
\begin{equation}\label{src:product-egf}
 C_{k,n-1}(x,y)
   =n![z^n w^{n-k}]\exp\bigl(wT(z)+\Ccal(z)\bigr).
\end{equation}
Here $C_{0,-1}=1$; for positive $k$ we set $C_{k,-1}=0$.
\end{lemma}

\begin{proof}
Put $\ell=n-k$.  Since $e^{\Ccal}=T'e^{-\sigma T}$, the
formal residue change $z=\phi(t)$ gives
\[
 \frac{n!}{\ell!}\operatorname*{Res}_{t=0}
      t^\ell e^{-\sigma t}\phi(t)^{-n-1}\,dt.
\]
We may work first with $\delta$ inverted.  Substitute
$u=e^{\delta t}-1$ and put $\beta=(n-1)x/\delta$.  The
display becomes
\[
 \frac{n!\delta^k}{\ell!}[u^n]
          (\log(1+u))^\ell(1+u)^\beta
 =\frac{\delta^k}{\ell!}
       \frac{d^\ell}{d\beta^\ell}\fall{\beta}{n}.
\]
Differentiating $\ell$ distinct linear factors and dividing by
$\ell!$ leaves the elementary symmetric polynomial of degree $k$
in the undifferentiated factors.  Thus the last expression is
\[
 e_k\bigl((n-1-i)x+iy:0\le i\le n-1\bigr)=C_{k,n-1}(x,y).
\]
Polynomiality extends the identity to $\delta=0$.
For $n=0$ it reduces to the empty-product convention.
\end{proof}

Extracting the factor $e^{wz}$ from \eqref{src:product-egf} yields
\begin{equation}\label{src:natural-source}
 W_{k,N}(x,y)=N![z^Nw^{N-k}]
       \exp\bigl(w(T(z)-z)+\Ccal(z)\bigr),
\end{equation}
where coefficients of negative powers of $w$ are zero.
Multiplication by $e^{wz}$ recovers
\[
 C_{k,n-1}(x,y)=\sum_{N=0}^{n}\binom nN W_{k,N}(x,y).
\]
This holds also when $n<k$, since both sides vanish.  Moreover,
$T-z$ has order at least two.  A term containing $w^{N-k}$ in
\eqref{src:natural-source} therefore has $z$-degree at least
$2(N-k)$, and hence
\begin{equation}\label{src:coarse-support}
 W_{k,N}=0\qquad\text{unless}\qquad k\le N\le2k.
\end{equation}
Consequently \eqref{src:Newton-bases} is a finite polynomial
expansion, of degree at most $2k$ in $d$.

\subsection{Positive differential equations and exact support}

\begin{theorem}\label{src:gamma-support}
For $k\ge1$ and $N\ge0$, there are nonnegative integers
$\Gamma_{k,p,N}$ such that
\begin{equation}\label{src:gamma-expansion}
 W_{k,N}(x,y)=\sum_{p=0}^{\lfloor k/2\rfloor}
          \Gamma_{k,p,N}\rho^p\sigma^{k-2p}.
\end{equation}
Their positive support is exactly
\begin{equation}\label{src:joint-support}
 \Gamma_{k,p,N}>0
 \quad\Longleftrightarrow\quad
 0\le p\le\lfloor k/2\rfloor,\quad
 k\le N\le2k-p,\quad p>0\ \hbox{or}\ N>k.
\end{equation}
The natural Schur coefficients are nonnegative integers, with
\begin{equation}\label{src:g-support}
 \begin{aligned}
  g_{k,0,N}>0&\quad\Longleftrightarrow\quad k+1\le N\le2k,\\
  g_{k,j,N}>0&\quad\Longleftrightarrow\quad k\le N\le2k
       &&(1\le j\le\lfloor k/2\rfloor).
 \end{aligned}
\end{equation}
In addition,
\begin{equation}\label{src:endpoints}
 \Gamma_{k,1,k}=(k-1)!(k-1)\quad(k\ge2),\qquad
 g_{k,j,2k}=(2k-1)!!
        \left(\binom kj-\binom k{j-1}\right).
\end{equation}
\end{theorem}

\begin{proof}
The equation $\phi''+\sigma\phi'+\rho\phi=0$, together with
$\phi(T)=z$, gives
\begin{equation}\label{src:ode}
 V'=\sigma V^2+\rho zV^3,\qquad V(0)=1;
 \qquad \Ccal'=\rho zV^2,\qquad\Ccal(0)=0.
\end{equation}
The recursions give nonnegative rational coefficients for $V$,
$T-z$, and $\Ccal$ in $\sigma,\rho,z$.  If
$\sigma$ and $\rho$ are assigned weights one and two, respectively,
then $[z^e]V$, $[z^{e+1}]T$, and $[z^e]\Ccal$ are homogeneous
of weight $e$.  Formula \eqref{src:natural-source} therefore gives
\eqref{src:gamma-expansion} with nonnegative rational coefficients.
Newton inversion expresses $W_{k,N}$ as an integral linear
combination of the polynomials $C_{k,n-1}$ for $0\le n\le N$.
Thus $W_{k,N}\in\mathbb Z[x,y]$ is symmetric, and
$\mathbb Z[x,y]^{\mathfrak S_2}=\mathbb Z[\sigma,\rho]$ makes
every $\Gamma_{k,p,N}$ an integer.

Induction in \eqref{src:ode} gives
\begin{equation}\label{src:V-support}
 [z^e\sigma^{e-2p}\rho^p]V>0
       \qquad(0\le p\le\lfloor e/2\rfloor).
\end{equation}
Indeed, if $e-2p>0$, the term $\sigma V^2$ uses the indicated
coefficient of $[z^{e-1}]V$; if $e=2p>0$, the term
$\rho zV^3$ uses the coefficient with indices $e-2,p-1$.
The induction starts at $V(0)=1$.  Thus weight $e$ and gamma
index $p$ occur in $T-z$ precisely when $e\ge\max(1,2p)$,
and in $\Ccal$ precisely when $p\ge1$ and $e\ge2p$.

Let $m=N-k$ be the number of factors from $T-z$ in a term of
\eqref{src:natural-source}.  Each such factor, with gamma index
$p_i$, has weight at least $p_i+1$; a factor from $\Ccal$ has
weight at least its gamma index.  Thus $k\ge m+p$, and
homogeneity also gives $k\ge2p$.  If $p=0$ and $k>0$, a factor
from $T-z$ is necessary, so $m>0$.  This proves necessity in \eqref{src:joint-support}.

Conversely, suppose $k\ge2p$, $0\le m\le k-p$, and $m+p>0$.
When $m=0$, a single term of $\Ccal$ realizes the required indices.
When $m\ge p$ and $m>0$, take $p$ factors from $T-z$ of gamma
index one and weight two, and $m-p$ factors of gamma index zero
and weight one.  Their total weight is $m+p\le k$.  When
$0<m<p$, distribute $p$ as a sum of $m$ positive indices and
use the corresponding factors of weights twice those indices;
the total weight is $2p\le k$.  In either case any excess weight
can be assigned to one factor.  The coefficients used in this construction are positive by
\eqref{src:V-support}, proving sufficiency.

For the Schur expansion, the elementary identity
\begin{equation}\label{src:ballot}
 \rho^p\sigma^{k-2p}
 =\sum_{j=p}^{\lfloor k/2\rfloor}
   \left(\binom{k-2p}{j-p}-\binom{k-2p}{j-p-1}\right)
            s_{(k-j,j)}(x,y)
\end{equation}
follows by taking successive differences of monomial coefficients.
All displayed coefficients are positive integers; binomial
coefficients with a negative lower index are zero.  The $p=0$
term of \eqref{src:gamma-expansion} gives every Schur index when
$k+1\le N\le2k$.  At $N=k$, the $p=1$ term gives exactly the
indices $j\ge1$.  This proves \eqref{src:g-support} and
integrality of the Schur coefficients.

Finally, setting $\rho=0$ in \eqref{src:ode} gives
$V=(1-\sigma z)^{-1}$ and $\Ccal=0$.  Hence
\[
 [\rho]\Ccal=[\rho]e^{\Ccal}
       =\int_0^z\frac{u}{(1-\sigma u)^2}\,du
       =\sum_{e\ge2}\frac{e-1}{e}\sigma^{e-2}z^e.
\]
This proves the first formula in \eqref{src:endpoints}.  At $N=2k$,
the power $w^k$ forces $k$ copies of the leading term
$\sigma z^2/2$ of $T-z$.  Thus
$W_{k,2k}=(2k)!\sigma^k/(2^k k!)$, and \eqref{src:ballot}
proves the second formula.
\end{proof}

Pascal's identity now gives
\begin{equation}\label{src:Pascal}
 B_{k,j,N}=g_{k,j,N}+g_{k,j,N+1}.
\end{equation}
In particular the original coefficients are nonnegative integers
and have no internal zeros.  Their exact positive supports are
$[k,2k]$ for $j=0$ and $[k-1,2k]$ for $j\ge1$.

\begin{remark}\label{src:degree-zero}
The degree-zero case is $C_{0,d}=1$.  Its only nonzero coefficients
are $\Gamma_{0,0,0}=g_{0,0,0}=B_{0,0,0}=1$.
\end{remark}

\subsection{The shifted quotient and its positive coefficients}

Set
\begin{equation}\label{src:Q}
 Q(t)=\frac{t}{\phi(t)}
      =e^{\sigma t/2}\frac{\delta t/2}{\sinh(\delta t/2)},
 \qquad H(z)=\log\frac{T(z)}z=\log Q(T(z)).
\end{equation}
The expressions are interpreted by their formal limits at $\delta=0$.

\begin{proposition}\label{src:shifted}
The polynomial $A_{k,j}(d)$ is divisible by $\fall{d+1}{k}$.
Define
\[
 \begin{gathered}
 P_{k,j}(d)=\frac{A_{k,j}(d)}{\fall{d+1}{k}},\qquad
 F_{k,j}(s)=P_{k,j}(s+k-1),\\
 \F_k(s;x,y)=\sum_jF_{k,j}(s)s_{(k-j,j)}(x,y).
 \end{gathered}
\]
Then
\begin{align}
 \sum_jP_{k,j}(d)s_{(k-j,j)}(x,y)
    &=[t^k]e^{-\sigma t}Q(t)^{d+2},\label{src:quotient}\\
 \F_k(s;x,y)&=[z^k]\exp\bigl(\Ccal(z)+sH(z)\bigr).
       \label{src:shifted-source}
\end{align}
In particular $P_{k,j}$ and $F_{k,j}$ have degree at most $k$.
\end{proposition}

\begin{proof}
The residue in the proof of Lemma~\ref{src:product-source} equals
\[
 \fall nk[t^k]e^{-\sigma t}Q(t)^{n+1}
       \qquad(n\ge k).
\]
The coefficient on the right is polynomial in $n$ of degree at
most $k$.  Agreement at all integers $n\ge k$ gives
\eqref{src:quotient} and the asserted divisibility, with $d=n-1$.
The derivative form of Lagrange inversion is
\[
 [z^k]T'(z)f(T(z))=[t^k]f(t)Q(t)^{k+1}.
\]
Apply it to $f(t)=e^{-\sigma t}Q(t)^s$ and use
$e^{\Ccal}=T'e^{-\sigma T}$ and $e^H=Q(T)$ to obtain
\eqref{src:shifted-source}.
\end{proof}

To study $H$, introduce the generalized Bernoulli polynomials
\begin{equation}\label{src:Bernoulli-def}
 B_n^{(\alpha)}(w)=n![t^n]
       \left(\frac{t}{e^t-1}\right)^\alpha e^{wt},
\end{equation}
where the power is defined by the formal logarithm with constant
term zero.

\begin{lemma}\label{src:Bernoulli-integer}
For $n\ge1$,
\begin{equation}\label{src:Bernoulli-integral}
 B_n^{(n+1)}(w)=\prod_{i=1}^n(w-i),\qquad
 B_n^{(n)}(w)=\int_0^1\prod_{i=1}^n(w+u-i)\,du.
\end{equation}
The polynomial $B_n^{(n)}$ has simple roots $r_1<\cdots<r_n$
satisfying
$i-1<r_i<i$ and $r_i+r_{n+1-i}=n$.
\end{lemma}

\begin{proof}
The substitution $u=e^t-1$ in the residue for
$B_n^{(n+1)}(w)/n!$ gives $[u^n](1+u)^{w-1}$ and proves the
product formula.  Integration of the generating function in
\eqref{src:Bernoulli-def} over $w+u$, $0\le u\le1$, lowers
the order from $n+1$ to $n$, giving the integral formula.
At each integer $i=0,\ldots,n$, its integrand has strict sign
$(-1)^{n-i}$ for $0<u<1$.  The alternating signs give a root in each $(i-1,i)$.
Since the degree is $n$, these roots are simple and exhaust the roots.  Replacing $u$ by $1-u$ in the integral gives
$B_n^{(n)}(n-w)=(-1)^nB_n^{(n)}(w)$ and the reflection of roots.
These integer-order facts are classical; see
\cite[Theorem~1 and equation~(7)]{Stampach}.
\end{proof}

\begin{lemma}\label{src:H-positive}
For $H_n(x,y)=[z^n]H(z)$, $n\ge1$, one has
\begin{equation}\label{src:H-Bernoulli}
 H_n=\frac1n[t^n]Q(t)^n
     =\frac{\delta^n}{n\,n!}
          B_n^{(n)}\left(\frac{nx}{\delta}\right).
\end{equation}
The symmetric form $H_n$ has $n$ distinct negative roots in $x/y$,
and
\begin{equation}\label{src:H-gamma}
 [\sigma^{n-2p}\rho^p]H_n>0
             \qquad(0\le p\le\lfloor n/2\rfloor).
\end{equation}
\end{lemma}

\begin{proof}
Lagrange inversion applied to $H=\log Q(T)$ gives
\[
 H_n=\frac1n[t^{n-1}]\frac{Q'(t)}{Q(t)}Q(t)^n
     =\frac1n[t^n]Q(t)^n.
\]
Writing $Q(t)=e^{xt}\delta t/(e^{\delta t}-1)$ proves
\eqref{src:H-Bernoulli}.  With the roots from
Lemma~\ref{src:Bernoulli-integer}, it becomes
\[
 H_n(x,y)=\frac1{n\,n!}
          \prod_{i=1}^n\bigl((n-r_i)x+r_i y\bigr).
\]
The roots $-r_i/(n-r_i)$ are distinct and negative.
The paired roots $r,n-r$ contribute
\[
 ((n-r)x+ry)(rx+(n-r)y)
        =r(n-r)\sigma^2+(n-2r)^2\rho.
\]
Both coefficients are positive.  In odd degree the central root
contributes $(n/2)\sigma$.  Expansion proves \eqref{src:H-gamma}.
\end{proof}

\begin{theorem}\label{src:shifted-positive}
Write
\[
 a_{k,j,r}=F_{k,j}^{(r)}(0),\qquad
 F_{k,j}(s)=\sum_{r=0}^{k}\frac{a_{k,j,r}}{r!}s^r.
\]
For $k\ge1$ and $0\le j\le\lfloor k/2\rfloor$,
\begin{equation}\label{src:derivative-positive}
 a_{k,j,r}>0\ (1\le r\le k),\qquad
 a_{k,0,0}=0,\qquad a_{k,j,0}>0\ (j\ge1).
\end{equation}
More precisely, if
\[
 D_{k,r,p}=[s^r z^k\sigma^{k-2p}\rho^p]
                   e^{\Ccal+sH},
\]
then its exact positive support is
\begin{equation}\label{src:derivative-gamma-support}
 D_{k,r,p}>0\quad\Longleftrightarrow\quad
 \begin{cases}
  1\le r\le k,\quad
  0\le p\le\min(\lfloor k/2\rfloor,k-r),&r\ge1,\\
  1\le p\le\lfloor k/2\rfloor,&r=0.
 \end{cases}
\end{equation}
All coefficients outside this support are zero.
\end{theorem}

\begin{proof}
The homogeneous coefficient of $s^r$ is
\begin{equation}\label{src:derivative-source}
 [s^r]\F_k(s;x,y)=\frac1{r!}[z^k]e^{\Ccal}H^r.
\end{equation}
Equations \eqref{src:ode} and \eqref{src:H-gamma} give its
nonnegative gamma expansion.  A term of $H_n$ with gamma index
$p$ has $n\ge\max(1,2p)$, and hence $p\le n-1$.
Every nonconstant term of $e^{\Ccal}$ has gamma index at least
one and degree at least twice that index.  A contribution to
\eqref{src:derivative-source} therefore satisfies
$2p\le k$ and $p\le k-r$.

Conversely, for $r\ge1$, distribute $p$ among $r$ factors of
$H$.  If $p\le r$, use $p$ indices equal to one and the others
zero; the least total degree is $r+p$.  If $p>r$, use $r$
positive indices with sum $p$; the least total degree is $2p$.
Thus the least degree is $\max(2p,r+p)$.  Under the displayed
conditions it is at most $k$, and any extra degree may be assigned
to a single factor, preserving strict positivity by
\eqref{src:H-gamma}.  The constant term of $e^{\Ccal}$ suffices.
For $r=0$, the strict support of $\Ccal$ established in the proof
of Theorem~\ref{src:gamma-support} gives exactly
$1\le p\le\lfloor k/2\rfloor$.

For $r\ge1$, the positive term with $p=0$ contributes to every
Schur index in \eqref{src:ballot}.  For $r=0$, the term with
$p=1$ contributes to every $j\ge1$, whereas the one-row
coefficient is zero.  This proves \eqref{src:derivative-positive}.
\end{proof}

For later use, the first source coefficients and the highest
derivatives are
\begin{equation}\label{src:normalizations}
 H_1=\frac\sigma2,\qquad
 H_2=\frac{5\sigma^2+4\rho}{24},\qquad
 a_{k,j,k}=2^{-k}\left(\binom kj-\binom k{j-1}\right).
\end{equation}
At $\sigma=1$, put $b=\rho$, $E=e^{\Ccal}$, and $R=2H/z$.
Then $E_0=R_0=1$, $E_1=0$, and $R_1=(5+4b)/12$.
Their coefficient polynomials are nonnegative and satisfy
\begin{equation}\label{src:normalized-degrees}
 \deg_b[z^\ell]E\le\lfloor\ell/2\rfloor,\qquad
 \deg_b[z^\ell]R\le\lfloor(\ell+1)/2\rfloor\le\ell
       \quad(\ell\ge1).
\end{equation}
Finally, \eqref{src:derivative-positive} implies that every real
zero of $P_{k,j}$ is at most $k-1$.  Equality occurs only for
$j=0$, where it is a simple zero.  We next fix a positive shift and study the zeros in the binary variables.

\section{Bernoulli polynomials of real order and binary roots}\label{bern:section}

With the normalization \eqref{src:Bernoulli-def}, the polynomial
$B_n^{(\alpha)}(w)$ is monic of degree $n$ in $w$ and has coefficients
polynomial in $\alpha$.  For integral $\alpha\ge n$, real-rootedness
follows from the root theorem for Bernoulli polynomials of the second
kind and differentiation; see \cite[Theorem~1 and Corollary~2]{Stampach}.
We extend the argument to real order.

\begin{theorem}\label{bern:continuous-order}
Let $n\ge1$ and $\alpha\in\mathbb R$.
\begin{enumerate}
\item If $\alpha\ge n$, then $B_n^{(\alpha)}$ has $n$ distinct
real roots, all in $(0,\alpha)$.
\item If $\alpha>n+1$, all its roots lie in $(1,\alpha-1)$.
\item At $\alpha=n+1$,
\begin{equation}\label{bern:integer-product}
 B_n^{(n+1)}(w)=\prod_{i=1}^n(w-i).
\end{equation}
\end{enumerate}
\end{theorem}

\subsection{A cut integral and its signs}

On $\Omega=\mathbb C\setminus(-\infty,-1]$, set
\begin{equation}\label{bern:psi}
 \psi(u)=\frac{\log(1+u)}u
        =\int_0^1\frac{dv}{1+vu},\qquad \psi(0)=1,
\end{equation}
using the principal logarithm.  The integral shows that $\psi$
has negative imaginary part in the upper half-plane and positive
imaginary part in the lower half-plane; it is positive on the
real interval $(-1,\infty)$.  Thus $\psi$ has no zeros on
$\Omega$.  Its analytic logarithm is normalized by
$\log\psi(0)=0$, and $\psi^\beta$ is defined using this logarithm.

The residue substitution $u=e^t-1$ in
\eqref{src:Bernoulli-def} gives, for every integer $N\ge1$,
\begin{equation}\label{bern:residue}
 \frac{B_N^{(N+1+\beta)}(w)}{N!}
    =[u^N](1+u)^{w-1}\psi(u)^\beta.
\end{equation}
Indeed, after $dt=du/(1+u)$ is included, the factor involving
$t$ and $u$ in the residue is
$u^{-N-1}(t/u)^\beta$.  Hence \eqref{bern:residue} also holds formally in $w,\beta$.

For $t>1$, the upper boundary value on the cut satisfies
\begin{equation}\label{bern:cut-argument}
 \begin{aligned}
 \psi(-t+i0)&=\frac{\log(t-1)+i\pi}{-t},\\
 \arg\psi(-t+i0)
     &=\arg\bigl(\log(t-1)+i\pi\bigr)-\pi\in(-\pi,0).
 \end{aligned}
\end{equation}
It follows that
\begin{equation}\label{bern:imaginary-sign}
 \operatorname{Im}\psi(-t+i0)^\beta
 \begin{cases}
  <0,&0<\beta\le1,\\
  >0,&-1<\beta<0.
 \end{cases}
\end{equation}
The specified branch determines these strict signs.

\begin{lemma}\label{bern:cut-signs}
If $-1<\beta\le1$, $\beta\ne0$, and $0\le i\le N-1$ is an
integer, then
\begin{equation}\label{bern:cut-integral}
 \frac{B_N^{(N+1+\beta)}(i+1)}{N!}
 =\frac{(-1)^{N+i-1}}{\pi}
   \int_1^\infty\frac{(t-1)^i}{t^{N+1}}
          \operatorname{Im}\psi(-t+i0)^\beta\,dt.
\end{equation}
In particular,
\begin{equation}\label{bern:integer-signs}
 \operatorname{sgn}B_N^{(N+1+\beta)}(i+1)
 =\begin{cases}
  (-1)^{N+i},&0<\beta\le1,\\
  (-1)^{N+i-1},&-1<\beta<0.
 \end{cases}
\end{equation}
\end{lemma}

\begin{proof}
At $w=i+1$, the factor $(1+u)^{w-1}$ in
\eqref{bern:residue} is a polynomial, so it introduces no cut.
Deform the coefficient contour for
$(1+u)^i\psi(u)^\beta u^{-N-1}$ around $(-\infty,-1]$.
On a large circle of radius $R$, the absolute integral is
\[
 O\bigl(R^{i-\beta-N}(\log R)^\beta\bigr)=o(1),
\]
since $i-\beta-N\le-1-\beta<0$.  This estimate holds also
when $\beta<0$: uniformly on the circle,
$|\log(1+u)|$ is comparable to $\log R$.  On a circle of
radius $\varepsilon$ about $-1$, the corresponding bound is
\[
 O\bigl(\varepsilon^{i+1}|\log\varepsilon|^\beta\bigr)=o(1).
\]
The cut integral is absolutely convergent.  At infinity an
absolute majorant is
$O(t^{i-N-1-\beta}(\log t)^\beta)$, whose power of $t$ is
strictly less than $-1$; at $t=1$ it is
$O((t-1)^i|\log(t-1)|^\beta)$.

The lower boundary value is the conjugate of the upper one.
Their difference, including the factors
$(1-t)^i$ and $(-t)^{-N-1}$, gives
\eqref{bern:cut-integral}.  Inserting
\eqref{bern:imaginary-sign} proves \eqref{bern:integer-signs}.
\end{proof}

\subsection{Root intervals and the real-order theorem}

The generating function gives the reflection and differentiation
identities
\begin{equation}\label{bern:reflection-derivative}
 \begin{aligned}
 B_N^{(\alpha)}(\alpha-w)&=(-1)^N B_N^{(\alpha)}(w),\\
 \frac{d^{N-n}}{dw^{N-n}}B_N^{(\alpha)}(w)
     &=\frac{N!}{n!}B_n^{(\alpha)}(w)
          \qquad(1\le n\le N).
 \end{aligned}
\end{equation}

\begin{lemma}\label{bern:unit-bands}
Fix $N\ge1$ and write $\alpha=N+1+\beta$.
\begin{enumerate}
\item If $0<\beta\le1$, $B_N^{(\alpha)}$ has one simple
root in each interval $(i,i+1)$, $1\le i\le N-1$, and one
in $(N,N+\beta)$.  All roots lie in $(1,\alpha-1)$.
\item If $-1<\beta<0$, it has one simple root in
$(1+\beta,1)$ and one in each $(i,i+1)$,
$1\le i\le N-1$.  All roots lie in $(\alpha-N,N)$.
\item If $\beta=0$, its roots are $1,2,\ldots,N$.
If $\beta=-1$, it has one simple root in each $(i-1,i)$,
$1\le i\le N$.
\end{enumerate}
\end{lemma}

\begin{proof}
For $\beta>0$, the values at $1,\ldots,N$ alternate in sign
by Lemma~\ref{bern:cut-signs}, and the value at $N$ is negative.
Monicity gives one more sign change between $N$ and $+\infty$.
The resulting $N$ disjoint intervals contain all the roots,
which are therefore simple.  The smallest root is greater than
one, so reflection in \eqref{bern:reflection-derivative} places
the largest below $\alpha-1=N+\beta$.

For $\beta<0$, the value at one has sign $(-1)^{N-1}$,
opposite to the sign at $-\infty$.  The other $N-1$ sign
changes occur between consecutive integers from one to $N$.
These $N$ sign changes account for all roots and establish simplicity.
The largest is less than $N$, so reflection places the smallest
above $\alpha-N=1+\beta$.

The cases $\beta=0$ and $\beta=-1$ are the product and integral
formulas of Lemma~\ref{src:Bernoulli-integer}, respectively.
\end{proof}

\begin{proof}[Proof of Theorem~\ref{bern:continuous-order}]
First let $\alpha\ge n$ be nonintegral.  Set
$N=\lfloor\alpha\rfloor\ge n$ and
$\beta=\alpha-N-1\in(-1,0)$.  Lemma~\ref{bern:unit-bands}
gives $N$ simple roots of $B_N^{(\alpha)}$ in $(0,\alpha)$.
Apply $N-n$ argument derivatives and Rolle's theorem, using
\eqref{bern:reflection-derivative}.  Each derivative has one
simple root in every interval between consecutive roots of its
predecessor, so $B_n^{(\alpha)}$ has $n$ simple roots in the
same open interval.  If $\alpha=M$ is integral, start instead
with $B_M^{(M)}$, whose roots lie in $(0,M)$ by
Lemma~\ref{src:Bernoulli-integer}, and differentiate $M-n$
times.

For the sharper assertion, suppose $\alpha>n+1$ is nonintegral.
Choose $N=\lfloor\alpha\rfloor-1\ge n$ and
$\beta=\alpha-N-1\in(0,1)$.  The first part of
Lemma~\ref{bern:unit-bands} now starts the same argument with
all roots in $(1,\alpha-1)$.  If $\alpha=M$ is integral,
start with
$B_{M-1}^{(M)}(w)=\prod_{i=1}^{M-1}(w-i)$.
Since $M-1>n$, at least one derivative is taken, and Rolle's
theorem puts every resulting root strictly inside $(1,M-1)$.
The product at $\alpha=n+1$ was already proved in
Lemma~\ref{src:Bernoulli-integer}.
\end{proof}

\subsection{Factorization of the shifted source}

Substituting $Q(t)=e^{xt}\delta t/(e^{\delta t}-1)$ in
\eqref{src:quotient} gives
\begin{equation}\label{bern:source}
 \F_k(s;x,y)=\frac{(x-y)^k}{k!}
       B_k^{(k+1+s)}\left(1+(k-1+s)\frac{x}{x-y}\right).
\end{equation}
Both sides are polynomial in $x,y$, with the specialization $x=y$
understood by continuity.

\begin{corollary}\label{bern:binary-roots}
For every $k\ge1$ and every real $s>0$, the polynomial
$\F_k(s;x,1)$ has degree $k$ and $k$ distinct negative roots.
Equivalently, $p\mapsto\F_k(s;p,p-1)$ has $k$ distinct roots
in $(0,1)$.  The gamma polynomial $G_{k,s}$ defined by
\begin{equation}\label{bern:gamma-polynomial}
 \F_k(s;x,y)=\sigma^kG_{k,s}(\rho/\sigma^2)
\end{equation}
has degree $\lfloor k/2\rfloor$ and distinct negative roots.
At $s=0$, for $k\ge2$,
\begin{equation}\label{bern:zero-shift}
 \F_k(0;x,y)=\frac1{k!}
       \prod_{i=0}^{k-1}\bigl(ix+(k-1-i)y\bigr).
\end{equation}
For $k=1$, $\F_1(s;x,y)=s(x+y)/2$.
\end{corollary}

\begin{proof}
Put $d=k-1+s>0$ and $\alpha=k+1+s=d+2$.  By
Theorem~\ref{bern:continuous-order}, the Bernoulli roots
$w_1,\ldots,w_k$ are distinct and belong to $(1,d+1)$.
Monicity and \eqref{bern:source} give
\begin{equation}\label{bern:positive-factorization}
 \F_k(s;x,y)=\frac1{k!}\prod_{i=1}^k
       \bigl((d+1-w_i)x+(w_i-1)y\bigr).
\end{equation}
Both coefficients in every factor are positive.  Its root in
$x/y$ is $-(w_i-1)/(d+1-w_i)$, and these roots are distinct.
Equivalently, $p=(w_i-1)/d$ belongs to $(0,1)$.

Symmetry pairs the negative binary roots under inversion.
For odd $k$ the unpaired root is $-1$ and contributes a factor
$\sigma$; for even $k$ there is no unpaired root.  A reciprocal
pair contributes a positive multiple of
$a\sigma^2+b\rho$, with $a,b>0$.  Each pair therefore gives
a negative root of $G_{k,s}$.  Distinct pairs give distinct
roots because $t\mapsto t/(1+t)^2$ is one-to-one on $(-1,0)$.
This proves the gamma assertion, including its degree.

At $s=0$, insert \eqref{bern:integer-product} in
\eqref{bern:source} to obtain \eqref{bern:zero-shift}.
The formula for $k=1$ follows from $H_1=\sigma/2$ and
$[z]\Ccal=0$.
\end{proof}

\subsection{Ultra-log-concavity in the Schur index}

\begin{corollary}\label{bern:Schur-ULC}
Let $h=\lfloor k/2\rfloor$.  For every real $s>0$,
$F_{k,j}(s)>0$ for $0\le j\le h$, and
\begin{equation}\label{bern:Schur-inequality}
 F_{k,j}(s)^2>
   \frac{(j+1)(k+2-j)}{j(k+1-j)}
       F_{k,j-1}(s)F_{k,j+1}(s)
       \qquad(1\le j<h).
\end{equation}
Thus $F_{k,j}(s)/\binom{k+1}{j}$ is strictly log-concave
throughout the interior of its support.  At $s=0$, the
corresponding weak inequalities hold.
\end{corollary}

\begin{proof}
Positivity follows from Theorem~\ref{src:shifted-positive}.
Write $C(t)=\F_k(s;t,1)=\sum_{i=0}^k c_it^i$ and set
$D(t)=(1-t)C(t)=\sum_{i=0}^{M}d_it^i$, where $M=k+1$.
Since $s_{(k-j,j)}(t,1)=t^j+\cdots+t^{k-j}$, Schur
extraction gives
\begin{equation}\label{bern:Schur-difference}
 d_j=c_j-c_{j-1}=F_{k,j}(s)\qquad(0\le j\le h),
 \qquad c_{-1}=0.
\end{equation}
By Corollary~\ref{bern:binary-roots}, $D$ has $M$ distinct
real roots: the $k$ negative roots of $C$ and the root one.

Newton's inequalities apply here even though the later coefficients
of $D$ change sign.  The following argument also gives strictness.  Fix $1\le j<h$ and put
$m=M-j+1$.  Differentiate $D$ exactly $j-1$ times, reverse the
resulting polynomial of degree $m$, and differentiate $m-2$
more times.  Since $d_{j-1}>0$, the first polynomial has
nonzero constant term, so reversal preserves its degree and its
simple real roots.  Rolle's theorem then gives a quadratic with
two distinct real roots.  Its coefficients, from the quadratic
term to the constant term, are
\[
 \frac{m!}{2}(j-1)!d_{j-1},\qquad
 (m-1)!j!d_j,\qquad
 \frac{(m-2)!}{2}(j+1)!d_{j+1}.
\]
The discriminant is strictly positive.  Dividing by the positive
factorials yields
\[
 d_j^2>
 \frac{(j+1)(M-j+1)}{j(M-j)}d_{j-1}d_{j+1}.
\]
Use \eqref{bern:Schur-difference} to obtain
\eqref{bern:Schur-inequality}; the weak inequalities at $s=0$
follow by continuity.
\end{proof}

Inequality \eqref{bern:Schur-inequality} concerns the Schur index
at fixed shift.  For the Newton index, we use the derivative source
\eqref{src:derivative-source} and the transfer results below.

\section{From derivatives to the original Newton coefficients}
\label{sec:bridge}

Fix the Schur index.  We transfer log-concavity in the derivative
order to log-concavity in the Newton degree.  The underlying
convolution principle is classical; the factorial weights and the
change between the two bases in \eqref{intro:newton} yield an explicit
strengthening of the Newton inequalities.

\subsection{A convolution-array principle}

A nonnegative sequence with interval support is log-concave precisely
when every minor of order two of its Toeplitz matrix is nonnegative.
This follows by comparing successive ratios, with interval support
settling the cases involving zeros.  The Cauchy--Binet formula then
shows that convolution preserves log-concavity.  Reflected sequences
are allowed, as are infinite sequences when the relevant sums are finite.

\begin{lemma}[Transposed convolution powers]\label{bridge:powers}
Let $(a_n)_{n\geq0}$ be finite, nonnegative and log-concave, with interval
support.  Let $q(t)=\sum_{n\geq1}q_nt^n$, where $(q_n)_{n\geq0}$,
with $q_0=0$, has the same properties apart from possible infinitude.
Then
\[
 b_m=\sum_{n\geq0}a_n[t^n]q(t)^m\qquad(m\geq0)
\]
is nonnegative, log-concave and interval-supported.
\end{lemma}
\begin{proof}
Fix $m\geq1$ and set
\[
 A_i=\sum_{n\geq0}a_{i+n}[t^n]q(t)^{m-1}\qquad(i\in\mathbb Z).
\]
This is a convolution of $a$ with a reflected convolution power, so
$A$ is log-concave.  Decreasing ratios imply
$A_iA_j\geq A_0A_{i+j}$ for $i,j\geq0$; this is immediate if
$A_0=0$.  Consequently
\[
 b_m^2=\Bigl(\sum_iq_iA_i\Bigr)^2
 \geq A_0\sum_{i,j}q_iq_jA_{i+j}=b_{m-1}b_{m+1}.
\]
If the nonempty supports of $a$ and $q$ are $[A,B]$ and $[c,d]$,
respectively, where $c\geq1$ and $d$ may be infinite, then $b_m>0$
for $m\geq1$ exactly when $mc\leq B$ and $md\geq A$.  These are
interval conditions.  If $b_0=a_0>0$, any nonempty positive support
at $m\geq1$ starts at $1$.  The empty cases are immediate.
\end{proof}

Lemma~\ref{bridge:powers} is a form of the convolution-array theorem
of F\"oldes--Major \cite[Proposition~3]{FM}.  More explicitly, if
$M=\max\operatorname{supp}a$, $c_h=a_{M-h}$ and $r(t)=q(t)/t$, then
\[
 b_m=[t^{M-m}]\Bigl(\sum_hc_ht^h\Bigr)r(t)^m.
\]
Thus $b$ is a reversed row of their convolution array.

\begin{corollary}\label{bridge:derivative}
If the derivative sequence $(F^{(r)}(0))_{r\geq0}$ of a real polynomial
is nonnegative, log-concave and interval-supported, then its forward
difference sequence $(\Delta^mF(0))_{m\geq0}$ has these properties.
\end{corollary}
\begin{proof}
Apply Lemma~\ref{bridge:powers} with $q(t)=e^t-1$.  Its nonzero
coefficients $1/n!$ form a log-concave sequence.  The identity
\[
 [t^r](e^t-1)^m=\frac{m!}{r!}
 \genfrac{\{} {\}}{0pt}{}{r}{m}
\]
identifies the output with $\Delta^mF(0)$.
\end{proof}

\subsection{An exact factor for the adjacent sum}

\begin{lemma}\label{bridge:weighted}
Let $k\geq1$ and let $b$ be a finite nonnegative log-concave sequence
with interval support, extended by zero.  Set
\[
 \alpha_m=\frac{m}{k+m},\qquad
 D_m=b_m+\alpha_mb_{m-1}\quad(m\geq0).
\]
Then $D$ is log-concave and interval-supported.  Its inequality at
$m\geq1$ is strict whenever $b_{m-1}b_m>0$.
\end{lemma}
\begin{proof}
For $m\geq1$, direct expansion gives
\begin{align*}
 D_m^2-D_{m-1}D_{m+1}
 ={}&(b_m^2-b_{m-1}b_{m+1})\\
 &+\alpha_{m-1}\alpha_{m+1}(b_{m-1}^2-b_{m-2}b_m)\\
 &+(\alpha_m^2-\alpha_{m-1}\alpha_{m+1})b_{m-1}^2\\
 &+\alpha_{m-1}(b_{m-1}b_m-b_{m-2}b_{m+1})\\
 &+(2\alpha_m-\alpha_{m-1}-\alpha_{m+1})b_{m-1}b_m.
\end{align*}
The three bracketed differences are nonnegative by log-concavity.
The remaining coefficients are positive, since
\begin{align*}
 \alpha_m^2-\alpha_{m-1}\alpha_{m+1}
 &=\frac{k(k+2m)}{(k+m)^2((k+m)^2-1)},\\
 2\alpha_m-\alpha_{m-1}-\alpha_{m+1}
 &=\frac{2k}{(k+m-1)(k+m)(k+m+1)}.
\end{align*}
The assertion about support follows directly from the definition.
\end{proof}

\begin{theorem}[Full-row transfer]\label{bridge:full}
Let $k\geq1$ and
\[
 A(d)=\fall{d+1}{k}F(d-k+1)
     =\sum_Ng_N\binom{d+1}{N}=\sum_NB_N\binom dN.
\]
If $(F^{(r)}(0))_r$ is nonnegative, log-concave and interval-supported,
then both $g$ and $B$ have these properties.  Moreover, for $N\geq k$,
\begin{equation}\label{bridge:quantitative}
 B_N^2\geq
 \left(1+\frac{k}{(N-k+1)(N+2)}\right)B_{N-1}B_{N+1}.
\end{equation}
For the Chern quotient $F=F_{k,j}$, this inequality is strict at every
$k\leq N\leq2k-1$ whenever the derivative hypothesis holds.
\end{theorem}
\begin{proof}
Put $b_m=\Delta^mF(0)$.  Corollary~\ref{bridge:derivative} applies.
Expanding $F$ in its Newton basis and writing $n=d+1$ gives
\[
 A(n-1)=\fall n k\sum_mb_m\binom{n-k}{m}
       =\sum_m\frac{(k+m)!}{m!}b_m\binom n{k+m}.
\]
Thus $g_{k+m}=w_mb_m$, where $w_m=(k+m)!/m!$.
Pascal's identity gives $B_N=g_N+g_{N+1}$, and hence
\[
 B_{k-1+m}=w_m(b_m+\alpha_mb_{m-1})=w_mD_m.
\]
Now use Lemma~\ref{bridge:weighted} and
\[
 \frac{w_m^2}{w_{m-1}w_{m+1}}
 =\frac{(k+m)(m+1)}{m(k+m+1)}
 =1+\frac{k}{m(k+m+1)}.
\]
The same ratio proves log-concavity of $g$.

For the Chern source, Theorem~\ref{src:shifted-positive} gives
$a_r=F_{k,j}^{(r)}(0)>0$ for $1\leq r\leq k$.  Also $a_0>0$ for
$j\geq1$, and $a_0=0$ for $j=0$.  The Stirling conversion therefore
gives the same statements for $b_0,\ldots,b_k$.  The strict part of
Lemma~\ref{bridge:weighted} applies at all the asserted centers except
possibly $j=0,N=k$.  At that center $B_{k-1}=0<B_k$, so the inequality
is again strict.
\end{proof}

\subsection{Transfer of terminal inequalities}

\begin{theorem}[Terminal transfer]\label{bridge:local}
Let $F(s)=\sum_{r=0}^ka_rs^r/r!$ with $a_r\geq0$, and define $g,B$
as in Theorem~\ref{bridge:full}.  Let $0\leq Q\leq k-1$ and
$L=k-Q-1$.  Suppose $(a_r)_{r=L}^k$ is log-concave and
interval-supported.  Then \eqref{bridge:quantitative} holds at every
$N\geq2k-Q$.  It is strict whenever $g_N>0$.
\end{theorem}
\begin{proof}
Set $\widetilde a_r=a_r$ for $r\geq L$ and $0$ otherwise.  This
truncated sequence is log-concave; its first inequality is automatic.
Corollary~\ref{bridge:derivative} gives a log-concave sequence
$\widetilde b$.  Since $[t^r](e^t-1)^m=0$ for $r<m$, one has
$\widetilde b_m=b_m$ for $m\geq L$.  Therefore
$b_m^2\geq b_{m-1}b_{m+1}$ at every $m\geq k-Q$.
With $\eta_m=(k+m)(m+1)/(m(k+m+1))$, this implies
\[
 g_{k+m}^2\geq\eta_mg_{k+m-1}g_{k+m+1}.
\]
Fix $m\geq k-Q\geq1$, and abbreviate
\[
 \begin{gathered}
 A=g_{k+m-1},\quad U=g_{k+m},\quad V=g_{k+m+1},\quad D=g_{k+m+2},\\
 \alpha=\eta_m,\qquad\beta=\eta_{m+1}.
 \end{gathered}
\]
The preceding inequalities give $U^2\geq\alpha AV$,
$V^2\geq\beta UD$, and $AD\leq UV/(\alpha\beta)$.
For $UV=0$, the last assertion follows from the interval support of
the truncated output on these four indices.  Thus
\[
 \beta(A+U)(V+D)
 \leq\frac\beta\alpha U^2+
       \left(\beta+\frac1\alpha\right)UV+V^2.
\]
Both coefficients in the resulting lower bound
\[
 (U+V)^2-\beta(A+U)(V+D)
 \geq\left(1-\frac\beta\alpha\right)U^2+
       \left(2-\beta-\frac1\alpha\right)UV
\]
are positive:
\begin{align*}
 1-\frac\beta\alpha
 &=\frac{k(k+2m+2)}{(k+m)(m+1)^2(k+m+2)},\\
 2-\beta-\frac1\alpha
 &=\frac{2k}{(k+m)(m+1)(k+m+2)}.
\end{align*}
Finally, with $N=k+m$,
$\beta=1+k/((N-k+1)(N+2))$.  Zero extension includes the upper
support boundary and proves the claim.
\end{proof}

Theorem~\ref{bridge:local} uses the derivative inequalities centered
at $a_{k-q}$, $1\leq q\leq Q$, and yields the Newton inequalities
starting at $N=2k-Q$.  Thus the adjacent sum preserves the first
terminal center.

\section{Uniform zero bounds near the balanced index}\label{near:section}\label{sec:near}

Write
\[
 k=2n+\epsilon,\qquad j=n-m,\qquad
 \omega=2m+\epsilon+1,
 \qquad \epsilon\in\{0,1\},\quad 0\le m\le n.
\]
The parameter $m$ measures distance from the balanced Schur index.
We obtain a zero bound for $P_{k,j}$ first when $n$ is at least
quadratic in $m$, and then in a range of linear width.

\begin{theorem}[An effective quadratic range]\label{near:quadratic}
If $m\ge2$ and $n\ge5(m+2)^2$, then every zero of
$P_{2n+\epsilon,n-m}$ belongs to the open disk
\[
                         |d|<\frac{7n}{5}.
\]
\end{theorem}

\begin{theorem}[A linear range]\label{near:linear}
If
\[
             n\ge10^{40},\qquad
             0<2m+\epsilon+1\le10^{-12}n,
\]
then every zero of $P_{2n+\epsilon,n-m}$ belongs to
$|d|<7n/5$.
In particular, this holds for both parities when
$0\le m\le n/(4\cdot10^{12})$ and $n\ge10^{40}$.
\end{theorem}

The proof of Theorem~\ref{near:linear} uses a two-variable contour
integral, with $d$ and the frequency allowed to be complex.
Appendix~\ref{near:linear-proof} gives the estimates, including the
limit at zero frequency.

\subsection{Character projection}\label{near:projection-section}

Proposition~\ref{src:shifted}, in the form
\eqref{src:quotient}--\eqref{src:Q}, gives
\begin{equation}\label{near:source}
 \sum_jP_{k,j}(d)s_{(k-j,j)}(x,y)
 =[t^k]e^{d(x+y)t/2}
 \left(\frac{(x-y)t/2}{\sinh((x-y)t/2)}\right)^{d+2}.
\end{equation}
The coefficient identity holds polynomially in $d$.  Set
\begin{equation}\label{near:H-def}
 \begin{gathered}
 H_h(d)=[u^{2h}](u/\sin u)^{d+2},\qquad
 J_h(v)=(2h+1)!v^{2h}H_h(1/v),\\
 \kappa_n=\frac1{4^n n!(n+1)!}.
 \end{gathered}
\end{equation}
The apparent singularity in the definition of $J_h$ is removable:
$H_h$ has degree $h$, so $J_h$ is a polynomial.
We use $(a)^{\overline\ell}=a(a+1)\cdots(a+\ell-1)$ for a rising
factorial, with value $1$ at $\ell=0$.

\begin{proposition}[Hahn projection]\label{near:hahn}
For $0\le m\le n$ one has
\begin{align}
 P_{2n+\epsilon,n-m}(d)
 &=\frac{\kappa_n c_{n,m}^{(\epsilon)}}{(n+2)^\epsilon}
   \sum_{h=0}^n\binom nh(2h+1)!
       d^{2n+\epsilon-2h}H_h(d)Q_m(h),\label{near:hahn-source}\\
 c_{n,m}^{(\epsilon)}
 &=\frac{2m+\epsilon+1}{\epsilon+1}
       \frac{\fall nm}{(n+\epsilon+2)^{\overline m}},
       \label{near:hahn-c}\\
 Q_m(h)
 &=\sum_{\ell=0}^{\min(m,h)}
 \frac{(-1)^\ell\fall m\ell
       (m+\epsilon+1)^{\overline\ell}\fall h\ell}
      {(3/2)^{\overline\ell}\fall n\ell\,\ell!}.
       \label{near:hahn-Q}
\end{align}
In particular $P_{k,j}$ has degree $k$ and positive leading coefficient.
\end{proposition}

\begin{proof}
Under $x=e^{i\theta}$ and $y=e^{-i\theta}$, the Schur polynomial
$s_{(k-j,j)}$ becomes the $SU(2)$ character
$U_{k-2j}(\cos\theta)$, where
\[
 U_q(v)=\sum_{a=0}^{\lfloor q/2\rfloor}
                 (-1)^a\binom{q-a}{a}(2v)^{q-2a}.
\]
Character orthogonality therefore extracts $P_{k,j}$ by integration
against $(2/\pi)\sin^2\theta\,U_{2m+\epsilon}(\cos\theta)$ on
$[0,\pi]$.  Expanding the two factors of the source reduces every term
to the beta integral
\begin{equation}\label{near:beta}
 \frac2\pi\int_0^\pi
       \cos^{2r}\theta\sin^{2h+2}\theta\,d\theta
 =\frac{(2r)!(2h+2)!}
        {2^{2(r+h)+1}r!(h+1)!(r+h+1)!}.
\end{equation}
For the balanced character $U_\epsilon=(2\cos\theta)^\epsilon$,
substitution of $r=n+\epsilon-h$ in \eqref{near:beta} gives the
coefficient
\[
           \frac{\kappa_n}{(n+2)^\epsilon}
                 \binom nh(2h+1)!.
\]

For the nonbalanced character, the formula for $U_q$ and the
binomial theorem give
\[
 \frac{U_{2m+\epsilon}(\sqrt{1-u})}
      {(2\sqrt{1-u})^\epsilon}
 =\frac{2m+\epsilon+1}{\epsilon+1}
  \sum_{\ell=0}^m
  \frac{(-1)^\ell\fall m\ell
        (m+\epsilon+1)^{\overline\ell}}
       {(3/2)^{\overline\ell}\ell!}\,u^\ell.
\]
In the beta integral the variable $u=\sin^2\theta$ has parameters
$h+3/2$ and $n+\epsilon-h+1/2$.  Thus the additional factor, apart
from $(2m+\epsilon+1)/(\epsilon+1)$, is the polynomial
\[
 E(h)=\sum_{\ell=0}^m
 \frac{(-1)^\ell\fall m\ell
       (m+\epsilon+1)^{\overline\ell}
       (h+3/2)^{\overline\ell}}
      {(3/2)^{\overline\ell}
       (n+\epsilon+2)^{\overline\ell}\ell!}.
\]
For $0\le r\le m$, taking $r$ forward differences in $h$ and using
$\Delta^r(h+a)^{\overline\ell}
=\fall\ell r(h+a+r)^{\overline{\ell-r}}$ yields
\[
 \Delta^r E(0)
 =\frac{(-1)^r\fall mr(m+\epsilon+1)^{\overline r}}
        {(3/2)^{\overline r}}
   \frac{(n-m+1)^{\overline{m-r}}}
        {(n+\epsilon+2)^{\overline m}}.
\]
Here the remaining sum is the terminating Vandermonde identity
\[
 \sum_{u=0}^{M}
 \frac{(-1)^u\fall Mu(b)^{\overline u}}
      {(c)^{\overline u}u!}
       =\frac{(c-b)^{\overline M}}{(c)^{\overline M}},
       \qquad c>b>0;
\]
it follows by integrating the binomial expansion of $(1-t)^M$
against the beta density with parameters $b,c-b$.
Since
$(n-m+1)^{\overline{m-r}}=\fall nm/\fall nr$,
Newton's finite expansion of the degree-$m$ polynomial $E$ gives
\[
 E(h)=\frac{\fall nm}{(n+\epsilon+2)^{\overline m}}Q_m(h).
\]
This proves \eqref{near:hahn-source}.  Its $h=0$ term is the unique
term of degree $k$ in $d$; all the factors in
\eqref{near:hahn-c} are positive.
\end{proof}

For $n\ge1$ define
\begin{equation}\label{near:reciprocal}
 \mathcal S_{n,m}^{(\epsilon)}(z)
    =\sum_{h=0}^n\binom nhJ_h(z/n)Q_m(h),\qquad
 A_n(z)=\sum_{h=0}^n\binom nhJ_h(z/n).
\end{equation}
Then $\mathcal S_{n,m}^{(\epsilon)}(0)=1$ and
\begin{equation}\label{near:reciprocal-source}
 P_{2n+\epsilon,n-m}(d)
 =\frac{\kappa_n c_{n,m}^{(\epsilon)}}{(n+2)^\epsilon}
   d^{2n+\epsilon}\mathcal S_{n,m}^{(\epsilon)}(n/d).
\end{equation}
Thus a zero-free disk for $\mathcal S$ gives a disk containing all
zeros of $P$, with a possible zero at $d=0$ already included.

\subsection{Coefficient majorants}

We write $V\preceq W$ when $W-V$ has nonnegative ordinary
coefficients.  Euler's sine product gives
\begin{equation}\label{near:H-recurrence}
 \begin{gathered}
 hH_h=(d+2)\sum_{q=1}^h\tau_qH_{h-q},\qquad
 \tau_q=\sum_{\nu\ge1}(\pi^2\nu^2)^{-q},\\
 \tau_1=\frac16,\qquad 0<\tau_q\le\frac{\tau_{q-1}}6\quad(q\ge2).
 \end{gathered}
\end{equation}
In particular, all coefficients of $H_h$ are positive.  Retaining
the $q=1$ term gives the lower estimate below.  For the upper estimate,
apply the recurrence for $H_{h-1}$ to the terms $q\ge2$:
\begin{equation}\label{near:H-successive}
 \frac{d+2}{6h}H_{h-1}\preceq H_h
       \preceq\frac{d+h+1}{6h}H_{h-1}.
\end{equation}
Induction consequently gives
\begin{equation}\label{near:H-majorant}
 \frac{d^h}{6^hh!}\preceq H_h(d)
       \preceq\frac{(d+2)^{\overline h}}{6^hh!}.
\end{equation}
Put
\[
 \begin{gathered}
 f(z)=(1-2z/3)^{-3/2},\qquad f(0)=1,\qquad
 b_h=(2/3)^h(3/2)^{\overline h},\\
 B_n(z)=\sum_{h=0}^n\binom nhb_h(z/n)^h.
 \end{gathered}
\]
Here $f$ denotes its analytic branch on $|z|<3/2$; the auxiliary
polynomial $B_n(z)$ in this subsection is unrelated to the Newton
coefficient $B_{k,j,N}$.

\begin{lemma}\label{near:comparison}
Let $0<r\le1$, set $R=r(1+2r/n)$, and assume $R<3/2$.
Then, throughout $|z|\le r$,
\begin{equation}\label{near:comparison-estimate}
 |A_n(z)-f(z)|
 \le\frac{2r^2f'(R)+r^2f''(r)}n.
\end{equation}
For $1\le\ell\le n$ one also has
\begin{equation}\label{near:shifted-majorant}
 \sum_{h=\ell}^n\binom{n-\ell}{h-\ell}J_h(r/n)
 \le b_\ell(r/n)^\ell
      \prod_{i=2}^{\ell+1}(1+ir/n)
      (1-2R/3)^{-\ell-3/2}.
\end{equation}
\end{lemma}

\begin{proof}
By \eqref{near:H-majorant}, the $h$th summand of $A_n(r)$ is at
most
\[
 \frac{b_hr^h}{h!}
 \prod_{i=0}^{h-1}
       (1-i/n)(1+(i+2)r/n).
\]
Every factor is nonnegative and
\begin{equation}\label{near:pairing}
 (1-i/n)(1+(i+2)r/n)
 =1+\frac{2r+(r-1)i}{n}-\frac{ri(i+2)}{n^2}
 \le1+2r/n.
\end{equation}
Extending the positive series proves $A_n(r)\le f(R)$.
Also, with $\fall nh/n^h$ set to zero for $h>n$,
\[
 0\le1-\frac{\fall nh}{n^h}\le\frac{h(h-1)}{2n}.
\]
For $h\le n$ this is the elementary product bound
$1-\prod_i(1-i/n)\le\sum_i i/n$; for $h>n$ the right-hand side
is at least $1$.  Therefore
\[
 0\le f(r)-B_n(r)\le\frac{r^2f''(r)}{2n}.
\]
Both $A_n-B_n$ and $f-B_n$ have nonnegative coefficients, so
\[
 |A_n(z)-f(z)|\le f(R)-f(r)+\frac{r^2f''(r)}n.
\]
Since $f'$ is increasing on $[r,R]$, this proves
\eqref{near:comparison-estimate}.

For \eqref{near:shifted-majorant}, set $h=\ell+u$ and use
\[
 J_h(r/n)\le b_h(r/n)^h
                  \prod_{i=2}^{h+1}(1+ir/n),\qquad
 b_{\ell+u}=b_\ell(2/3)^u(\ell+3/2)^{\overline u}.
\]
After taking out the factors with indices $2,\ldots,\ell+1$,
the remaining product pairs
$(1-(\ell+i)/n)$ with $(1+(\ell+i+2)r/n)$ for
$0\le i<u$.  Inequality~\eqref{near:pairing} bounds it by
$(1+2r/n)^u$.  Summing the resulting binomial series over all
$u\ge0$ proves \eqref{near:shifted-majorant}.
\end{proof}

\subsection{Proof of the quadratic range}

The cancellation
\[
 \binom nh\frac{\fall h\ell}{\fall n\ell}
                  =\binom{n-\ell}{h-\ell}
\]
allows us to apply \eqref{near:shifted-majorant} term by term to
\eqref{near:hahn-Q}.  The paired factorials satisfy the sharper
inequality
\begin{equation}\label{near:hahn-pairing}
 \fall m\ell(m+\epsilon+1)^{\overline\ell}
 =\prod_{i=0}^{\ell-1}
     \{m(m+\epsilon+1)-i(i+\epsilon+1)\}
 \le\{m(m+\epsilon+1)\}^{\ell},
\end{equation}
whose factors are nonnegative for $0\le\ell\le m$.
Consequently, on $|z|\le r$,
\begin{equation}\label{near:hahn-error}
 |\mathcal S_{n,m}^{(\epsilon)}(z)-A_n(z)|
 \le f(R)\sum_{\ell=1}^m\frac{\lambda^\ell}{\ell!}
 \le f(R)\frac{\lambda}{1-\lambda/2},
\end{equation}
provided $0\le\lambda<2$, where
\begin{equation}\label{near:lambda}
 \lambda=\frac{m(m+\epsilon+1)}n
       \frac{2r}{3(1-2R/3)}
       \left(1+\frac{(m+1)r}{n}\right).
\end{equation}
The last estimate uses $\ell!\ge2^{\ell-1}$ for $\ell\ge1$.
The estimate uses the absolute coefficients of the Hahn expansion;
no sign condition on $Q_m(h)$ is required.

\begin{proof}[Proof of Theorem~\ref{near:quadratic}]
Take $r=5/7$.  For fixed $m$, the bounds in
\eqref{near:comparison-estimate}, \eqref{near:hahn-error}, and
\eqref{near:lambda} decrease with $n$.  It is enough to evaluate
them at $n_0=5(m+2)^2$.  Put
\[
 q=\frac1{m+2}\le\frac14,\qquad
 h=1-\frac{20q^2}{77},\qquad p_0=1-2r/3=\frac{11}{21}.
\]
At $n=n_0$, one has $1-2R/3=p_0h$.  Using $\epsilon\le1$ gives
\begin{equation}\label{near:lambda-rational}
 \lambda\le\lambda_*
 :=\frac2{11}\frac{(1-2q)(1+q(1-q)/7)}h
 \le\frac2{11}\left(1-\frac{9q}{5}\right).
\end{equation}
The second inequality is the polynomial identity
\[
 h(1-9q/5)-(1-2q)(1+q(1-q)/7)
       =\frac{2q}{35}+\frac{13q^2}{77}+\frac{2q^3}{11}\ge0.
\]
In particular $0\le\lambda_*<2/11$.

Divide the sum of the two error estimates by $f(r)$.
Since $h^{-5/2}\le h^{-3}$ and $h^{-3/2}\le h^{-2}$, the result
is bounded by
\[
 E(q)=q^2\left(\frac{30}{77h^3}+\frac{75}{121}\right)
        +h^{-2}\frac{\lambda_*}{1-\lambda_*/2}.
\]
For $0\le q\le1/4$, the following are rational inequalities:
\[
 h\ge\frac{303}{308},\qquad h^{-3}<\frac{16}{15},\qquad
 h^{-2}\le1+\frac{11q^2}{20},\qquad
 \frac{30}{77}\frac{16}{15}+\frac{75}{121}<\frac{21}{20}.
\]
For example,
$h^{-2}-1\le(40/77)(308/303)^2q^2<(11/20)q^2$;
the other comparisons follow by multiplication by their positive
integer denominators.  The function $\phi(x)=x/(1-x/2)$ is
increasing and has $\phi'(x)\ge1$ on $[0,2)$.  As
$\phi(2/11)=1/5$, \eqref{near:lambda-rational} implies
\[
 \phi(\lambda_*)\le\frac15-\frac{18q}{55}.
\]
Multiplying by $1+11q^2/20$ and discarding the negative cubic term
now gives
\[
 E(q)\le\frac15-\frac{18q}{55}+\frac{29q^2}{25}
       \le\frac15-\frac{41q}{1100}<\frac15.
\]
Here the second inequality is $q^2\le q/4$.
Thus $|\mathcal S_{n,m}^{(\epsilon)}-f|<f(r)/5$ on the closed disk.
On its boundary,
\[
 |f(z)|\ge\left(\frac{21}{31}\right)^{3/2}
       >\frac15\left(\frac{21}{11}\right)^{3/2}
       =\frac{f(r)}5,
\]
where the strict inequality is $31^3<25\cdot11^3$.
Rouch\'e's theorem, including the strict boundary estimate, proves
that $\mathcal S_{n,m}^{(\epsilon)}$ is zero-free on $|z|\le5/7$.
Equation~\eqref{near:reciprocal-source} proves the assertion.
\end{proof}

\subsection{From zero disks to derivative inequalities}

\begin{lemma}[A derivative sector criterion]\label{near:sector}
Let $G$ be a real polynomial with positive leading coefficient and
all zeros in
\[
 K=\{z:\Re z\le0,\ |\Im z|\le-\Re z\}.
\]
Then $(G^{(r)}(0))_r$ is nonnegative, log-concave, and
interval-supported.  If all zeros lie in the interior of $K$, all
derivatives through the degree are positive, and every inner
log-concavity inequality is strict.  The same conclusions hold at
every real argument $t\ge0$.
\end{lemma}

\begin{proof}
Conjugate root pairs give factors
$s^2+2us+u^2+v^2$, with $u>0$ and $|v|\le u$; real zeros give
factors $s+u$, with $u\ge0$.  This proves the sign and support
assertions.  The sector is convex, so Gauss--Lucas keeps the zeros
of every nonconstant derivative in $K$.  If $V=G^{(r-1)}$ has no
zero at the origin, logarithmic differentiation gives
\begin{equation}\label{near:log-derivative}
 (V'(0))^2-V(0)V''(0)
       =V(0)^2\sum_{\zeta:V(\zeta)=0}\frac1{\zeta^2},
\end{equation}
with zeros counted with multiplicity.  Each real-zero contribution
is positive; a conjugate pair contributes
$2\Re(\zeta^{-2})\ge0$.  In the open sector every such
contribution is positive, and Gauss--Lucas places derivative zeros
in the convex hull of a finite set strictly inside that sector.
This proves strictness.  Zeros at the origin and boundary zeros are
handled by $G(s+\delta)$ followed by $\delta\downarrow0$.
Finally, replacing $s$ by $s+t$ moves every zero farther to the left.
\end{proof}

\begin{corollary}\label{near:band-consequences}
In either Theorem~\ref{near:quadratic} or
Theorem~\ref{near:linear}, the derivative sequence
$(F_{k,j}^{(r)}(t))_{r=0}^k$ is positive and strictly log-concave
at every inner index for every real $t\ge0$.  The ordinary Newton
row $(B_{k,j,N})_N$ has support $[k-1,2k]_{\mathbb Z}$ and is
log-concave; it also satisfies the quantitative inequalities of
Theorem~\ref{bridge:full}.
\end{corollary}

\begin{proof}
The zeros of $F_{k,j}$ lie in the disk with center $-(k-1)$ and
radius $7n/5$.  For $n\ge50$,
\[
 (k-1)^2-2(7n/5)^2
 \ge(2n-1)^2-98n^2/25
 =\frac{2}{25}n(n-50)+1>0.
\]
Hence the disk is strictly inside $K$.  Both hypotheses in the
corollary imply $n\ge50$, and the leading coefficient is positive
by Proposition~\ref{near:hahn}.  Lemma~\ref{near:sector} and
Theorem~\ref{bridge:full} apply.  All derivatives through degree
$k$ are positive, so all finite differences of orders $0,\ldots,k$
are positive; the support assertion follows from the exact transfer
formula in that theorem.
\end{proof}

\subsection{The balanced family}

For the balanced index, a smaller zero disk reduces the proof to four
finite degrees.

\begin{theorem}\label{near:balanced}
For $k\ge2$ the balanced derivative row
$(F_{k,\lfloor k/2\rfloor}^{(r)}(0))_{r=0}^k$ is positive and
strictly log-concave at every inner index.  For $k=1$ the row is
$(0,1/2)$.  At every real $t\ge0$ the derivative row remains
nonnegative and log-concave with interval support, and all its
inner inequalities for $k\ge2$ remain strict.
\end{theorem}

\begin{proof}
Put
\begin{equation}\label{near:S-def}
 S_n(d)=\sum_{h=0}^n\binom nh(2h+1)!d^{2n-2h}H_h(d).
\end{equation}
Proposition~\ref{near:hahn} gives
\begin{equation}\label{near:balanced-formulas}
 P_{2n,n}=\kappa_nS_n,\qquad
 P_{2n+1,n}=\frac{\kappa_nd}{n+2}S_n.
\end{equation}
The two estimates in \eqref{near:H-successive}, applied separately
to $(h+1)H_h$ and $dH_h$, show that
\begin{equation}\label{near:balanced-coefficient}
 (h+1-d)H_h+\frac{d^2}{4h+2}H_{h-1}\succeq0
       \qquad(h\ge1).
\end{equation}
Indeed a coefficientwise lower bound is $H_{h-1}$ times
\[
 \frac{h+1}{3h}
     +\left(\frac1{4h+2}-\frac1{6h}\right)d^2,
\]
which has nonnegative coefficients.

Let $R_n(z)=\sum_{h=0}^n\binom nhJ_h(z)$, the reversal of $S_n$.
Reversing \eqref{near:balanced-coefficient} and using
$h\binom nh=(n-h+1)\binom n{h-1}$ gives
\[
 ((n+1)z-1)R_n(z)
 =(z-1)+\sum_{h=1}^n\binom nh
       \{((h+1)z-1)J_h(z)+hzJ_{h-1}(z)\}.
\]
Every nonconstant coefficient of this expression is nonnegative.
Moreover, all coefficients of $S_n$ are positive: the positive
degree intervals $[2n-2h,2n-h]$, $0\le h\le n$, cover
$[0,2n]$.  Writing $P_{2n,n}(d)=\sum_i p_i d^i$ therefore gives
$0<p_i\le(n+1)p_{i+1}$ for $0\le i<2n$.
It follows that
\begin{equation}\label{near:balanced-disk}
             P_{2n,n}(d)=0\quad\Longrightarrow\quad |d|\le n+1.
\end{equation}
To prove this bound, rescale the variable so that the coefficients
$c_0,\ldots,c_D$ are nondecreasing.  At a zero $z$ with $|z|>1$,
the identity
\[
 (1-z)\sum_{i=0}^Dc_iz^i
 =c_0+\sum_{i=1}^D(c_i-c_{i-1})z^i-c_Dz^{D+1}
\]
would imply
$c_D|z|^{D+1}\le c_0+\sum_i(c_i-c_{i-1})|z|^i
\le c_D|z|^D$, a contradiction.  This proves \eqref{near:balanced-disk}.

For $n\ge5$ the shifted disk lies strictly inside $K$, because
\[
 (2n-1)^2-2(n+1)^2=2n^2-8n-1>0.
\]
Lemma~\ref{near:sector} proves the even-degree assertion in that
range.  The remaining four derivative rows are, up to the indicated
positive factors,
\begin{align*}
 k=2:\quad&\tfrac18(4,3,2),\\
 k=4:\quad&\tfrac1{288}(324,305,233,126,36),\\
 k=6:\quad&\tfrac1{41472}
       (165600,162692,136938,94389,49788,17820,3240),\\
 k=8:\quad&\tfrac1{4147200}
       (71794800,71527066,63200803,48437685,31262775,\\[-2pt]
       &\hspace{33mm}16294200,6423300,1701000,226800).
\end{align*}
These follow by substituting $n=1,2,3,4$ into
\eqref{near:H-recurrence} and \eqref{near:balanced-formulas}.
For the unscaled integer rows, the minimum inner Tur\'an difference
is respectively
$1$, $7488$, $156239280$, and $1436596560000$.
Thus every required finite inequality is strictly positive.

The odd-degree formula in \eqref{near:balanced-formulas} gives
\begin{equation}\label{near:balanced-odd}
 F_{2n+1,n}(s)=\frac{s+2n}{n+2}F_{2n,n}(s+1).
\end{equation}
Two preservation properties complete the argument.
A translation by $t\ge0$ sends a derivative row $a$ to
$b_r=\sum_{\ell\ge0}a_{r+\ell}t^\ell/\ell!$.
This is convolution with the reflection of the log-concave sequence
$(t^\ell/\ell!)_{\ell\ge0}$, so it preserves log-concavity and
interval support by the $2\times2$ Cauchy--Binet formula for
Toeplitz matrices.  Its adjacent minor includes the original
adjacent minor with coefficient $1$, which preserves strictness.
Multiplication by $\alpha+\beta s$, for $\alpha,\beta\ge0$,
sends $b_r$ to $e_r=\alpha b_r+\beta r b_{r-1}$, where entries
outside the support are zero.  Direct expansion gives
\begin{align*}
 e_r^2-e_{r-1}e_{r+1}
 ={}&\alpha^2(b_r^2-b_{r-1}b_{r+1})\\
 &+\alpha\beta(r-1)(b_rb_{r-1}-b_{r-2}b_{r+1})\\
 &+\beta^2\{r^2(b_{r-1}^2-b_{r-2}b_r)+b_{r-2}b_r\}\ge0.
\end{align*}
The middle term is nonnegative by decreasing consecutive ratios.
When $\alpha=2n>0$ and $\beta=1$, the last line is strictly
positive at every inner index of the extended positive row
(at $r=1$ it is $b_0^2$).  This proves the odd-degree assertion
for $n\ge1$, and $F_{1,0}(s)=s/2$ gives the remaining case.
The translation argument proves the assertion for every $t\ge0$.
\end{proof}

\subsection{The first nonbalanced family and finite closure}

For $m=1$, Proposition~\ref{near:hahn} reduces to
\begin{align}
 P_{2n,n-1}(d)
 &=\frac{\kappa_n n}{n+2}
                  \{4d^2S_{n-1}(d)-S_n(d)\},\label{near:first-even}\\
 P_{2n+1,n-1}(d)
 &=\frac{2\kappa_n nd}{(n+2)(n+3)}
                  \{2d^2S_{n-1}(d)-S_n(d)\}.
                  \label{near:first-odd}
\end{align}
Here $S_0=1$.  We first bound the zeros for all sufficiently large
$n$; the remaining degrees will be treated in
Proposition~\ref{comp:central}.

\begin{proposition}\label{near:first-roots}
For $n\ge40$, every zero of $P_{2n+\epsilon,n-1}$ lies in
$|d|<7n/5$.  For $n\ge50$, all derivative rows at nonnegative
real shifts are positive and strictly log-concave.
\end{proposition}

\begin{proof}
For $a\in\{2,4\}$ set
$G_{n,a}(z)=aR_{n-1}(z/n)-R_n(z/n)$, where $R_n$ is the
unscaled reversal used above.  The reciprocal identity is
\[
 G_{n,a}(z)=(z/n)^{2n}
       \{ad^2S_{n-1}(d)-S_n(d)\},\qquad d=n/z.
\]
Put $c_h=b_h/h!$, and let $a_{n,h}=\fall nh/n^h$ for $h\le n$,
and $a_{n,h}=0$ for $h>n$.  Retaining only the leading term of
each $H_h$ gives a polynomial $G^0_{n,a}$ whose coefficient of
$c_hz^h$ is $a_{n,h}(a-1-ah/n)$.  For all $h\ge0$,
\[
 0\le(a-1)-a_{n,h}(a-1-ah/n)
 \le\frac{(a-1)h(h-1)}{2n}+\frac{ah}{n}.
\]
For $h\le n$, write the middle expression as
$(a-1)(1-a_{n,h})+a_{n,h}ah/n$; for $h>n$ the claimed upper
bound follows from $h(h-1)/(2n)\ge1$.
Thus, on $|z|\le r$,
\[
 |G^0_{n,a}(z)-(a-1)f(z)|
 \le\frac{(a-1)r^2f''(r)}{2n}+\frac{ar f'(r)}n.
\]
The discarded coefficients of $J_h$ are nonnegative and their
multiplier in $G_{n,a}$ is
\[
 a\binom{n-1}h-\binom nh
       =\binom nh(a-1-ah/n).
\]
Its last factor belongs to $[-1,a-1]$, and hence has absolute
value at most $a-1$.  The proof of Lemma~\ref{near:comparison}
therefore gives, for $n\ge N$,
\begin{equation}\label{near:first-error}
 \begin{split}
 |G_{n,a}(z)-(a-1)f(z)|
 &\le\frac{2(a-1)r^2}{n}f'(r(1+2r/N))\\
 &\quad+\frac{(a-1)r^2}{n}f''(r)+\frac{ar}{n}f'(r).
 \end{split}
\end{equation}
Take $r=5/7$ and $N=40$, so that $r(1+2r/N)=145/196$.
Since
$f'(z)=(1-2z/3)^{-5/2}$ and
$f''(z)=\tfrac53(1-2z/3)^{-7/2}$, one has
\[
 f'(r)\le f'(145/196)<\frac{11}{2},\qquad
 f''(r)<\frac{33}{2},\qquad
 \min_{|z|\le r}|f(z)|>\frac59.
\]
The strict comparisons reduce respectively to
\begin{align*}
 121\cdot149^5-4\cdot294^5&=100102720733>0,\\
 9801\cdot11^7-100\cdot21^7&=10884908871>0,\\
 81\cdot21^3-25\cdot31^3&=5366>0.
\end{align*}
It follows from \eqref{near:first-error} that
\[
 |G_{n,4}-3f|\le\frac{5665}{98n}
          \le\frac{5665}{3920}<\frac53<|3f|,
\]
and
\[
 |G_{n,2}-f|\le\frac{2145}{98n}
          \le\frac{2145}{3920}<\frac59<|f|.
\]
Rouch\'e's theorem proves zero-freeness in the closed reciprocal
disk, and \eqref{near:first-even}--\eqref{near:first-odd} give the
asserted root bound.  The sector calculation in the proof of
Corollary~\ref{near:band-consequences} applies for $n\ge50$.
\end{proof}

A related root-region multiplier criterion is stated by Tang
\cite[Appendix~B.2]{Tang}.  The sector argument above establishes
derivative log-concavity, beginning at $n=50$ for the first
nonbalanced family.

\begin{theorem}[Complete central families]\label{near:central-families}
Let $k\ge1$ and $0\le j\le\lfloor k/2\rfloor$.  If
$0\le k-2j\le9$, then $(B_{k,j,N})_N$ is nonnegative,
interval-supported, and log-concave.  Its support is
$[k,2k]_{\mathbb Z}$ for $j=0$ and $[k-1,2k]_{\mathbb Z}$
for $j\ge1$.

For $m\in\{0,3,4\}$ in $k=2n+\epsilon$, $j=n-m$, the full
derivative row is also nonnegative and strictly log-concave at
every inner index, and this remains true at every real shift
$t\ge0$.  At $t=0$ all entries are positive for $j\ge1$; for
$j=0$ only the initial entry is zero.
\end{theorem}

\begin{proof}
The case $m=0$ follows from Theorem~\ref{near:balanced} and the
transfer theorem.  Table~\ref{tab:central-ranges} gives the remaining
analytic and finite ranges; both parities are included.

\begin{table}[htbp]
\caption{Analytic and finite ranges for the central families.}
\label{tab:central-ranges}
\centering\small
\begin{tabular}{@{}cccc@{}}
\toprule
$m$ & Analytic range & Finite range & Finite rows\\
\midrule
1 & $n\ge50$ & $1\le n\le49$ & 98 Newton rows\\
2 & $n\ge80$ & $2\le n\le99$ & 196 Newton rows\\
3 & $n\ge125$ & $3\le n\le124$ & 244 derivative rows\\
4 & $n\ge180$ & $4\le n\le179$ & 352 derivative rows\\
\bottomrule
\end{tabular}
\end{table}

Proposition~\ref{near:first-roots} treats the first analytic range;
Theorem~\ref{near:quadratic} treats the other three.  The exact
finite certificates and their source reconstruction are given in
Proposition~\ref{comp:central}.  These ranges exhaust every $n\ge m$.
For $m=3,4$ the certified derivative inequalities, followed by
Theorem~\ref{bridge:full}, prove every ordinary Newton inequality.
Translation preserves their strict derivative inequalities by the
argument in the proof of Theorem~\ref{near:balanced}.
The support assertions follow from the certified derivative or
Newton supports and the exact transfer formulas.
\end{proof}

For $m=1,2$, the finite certificates establish the Newton inequalities.
The derivative conclusions in Theorem~\ref{near:central-families}
apply to $m=0,3,4$.

\section{Schur cumulants and a uniform terminal Newton band}
\label{tail:section}
\label{sec:tail}

Write $q=k-r$ for the defect from the highest derivative order.
We prove a log-concavity estimate uniform in the Schur index.

\begin{theorem}\label{tail:main}
Let $q\ge1$ and $k\ge 2^{38}(q+1)$ be integers.  For every
$0\le j\le\lfloor k/2\rfloor$,
\begin{equation}\label{tail:derivative-margin}
 \log\frac{a_{k,j,k-q}^{\,2}}
 {a_{k,j,k-q-1}a_{k,j,k-q+1}}
 \ge \frac{1}{2(q+1)}.
\end{equation}
Consequently, if $Q\ge1$ and $k\ge2^{38}(Q+1)$, then
\begin{equation}\label{tail:newton-margin}
 B_{k,j,N}^{\,2}>
 \left(1+\frac{k}{(N-k+1)(N+2)}\right)
 B_{k,j,N-1}B_{k,j,N+1}
\end{equation}
for $2k-Q\le N\le2k-1$.
\end{theorem}

The theorem covers a number of terminal inequalities proportional to
$k$.  Its proof combines two estimates: a degree bound for the
logarithmic corrections, obtained from connected cumulants, and a
nonvanishing disk of radius proportional to $1/q$, obtained from
uniform factorial bounds.  Together they control the discrete
curvature of the logarithm.

\subsection{The Schur functional and its normalization}

In the source of \eqref{src:shifted-source}, set $x+y=1$ and $b=xy$, and
write
\begin{equation}\label{tail:ER}
 E(z,b)=e^{\Ccal(z,b)},\qquad
 R(z,b)=\frac{2H(z,b)}{z},\qquad
 c(b)=\frac{5+4b}{12}.
\end{equation}
The source positivity proved above gives
\begin{equation}\label{tail:source-coefficients}
 \begin{gathered}
 E_0=R_0=1,\quad E_1=0,\quad R_1=c(b),\\
 E_\ell,R_\ell\in\mathbb Q_{\ge0}[b],\qquad
 \deg_b E_\ell,\deg_b R_\ell\le\ell\quad(\ell\ge1).
 \end{gathered}
\end{equation}
Here and below a subscript on a source series denotes its ordinary
$z$-coefficient.

For independent formal variables $u,t$, put $\nu=u/(1+u)$ and define a
unital linear functional on $\mathbb Q[b]$ by
\begin{equation}\label{tail:Lambda}
 \Lambda_{u,t}(b^p)=w_p(u,t):=
 \frac{\displaystyle\prod_{i=0}^{p-1}(t-i\nu+i^2\nu^2)}
 {\displaystyle\prod_{h=1}^{2p}(1-h\nu)},\qquad w_0=1.
\end{equation}
We extend $\Lambda_{u,t}$ coefficientwise whenever each coefficient of
the input is a polynomial in $b$.  All formal expansions in this
section are at $u=0$.  In particular, each coefficient of $w_p$ is a
polynomial in $t$.  The functional need not have a positive representing measure.
At $t=0$ it reduces to evaluation at zero:
\begin{equation}\label{tail:t-zero}
 \Lambda_{u,0}(f(b))=f(0).
\end{equation}

To identify the Schur functional, write
\[
 f_{k,j}=\binom{k}{j}-\binom{k}{j-1}>0,\qquad
 t=\frac{j(k+1-j)}{(k+1)^2},\qquad u=\frac1k.
\]
For $0\le p\le q<k/2$, specialization of \eqref{tail:Lambda} gives
\begin{equation}\label{tail:Schur-weight}
 w_p(1/k,t)=
 \frac{\fall{j}{p}\fall{k-j+1}{p}}{\fall{k}{2p}}
 =\frac{f_{k-2p,j-p}}{f_{k,j}},
\end{equation}
where $f_{k-2p,j-p}=0$ if $p>j$.  Indeed, the two factors in the
numerator corresponding to $i$ multiply to
$j(k+1-j)-i(k+1)+i^2$, while $\nu=1/(k+1)$; the last equality follows
from the binomial formula for $f_{k,j}$.  Equivalently,
\[
 [s_{(k-j,j)}](x+y)^{k-2p}(xy)^p=f_{k-2p,j-p}.
\]
Restoring the homogeneous degree after $x+y=1$, the derivative source
therefore yields
\begin{equation}\label{tail:physical-source}
 v_q:=\frac{2^ka_{k,j,k-q}}{f_{k,j}}
 =2^q\Lambda_{1/k,t}\bigl([z^q]E(z,b)R(z,b)^{k-q}\bigr).
\end{equation}
Only moments with index at most $q<k/2$ occur, so all denominators in
this specialization are nonzero.

For every integer $q\ge0$, define
\begin{align}
 M_q(u,b)&=q!u^q[z^q]E(z,b)R(z,b)^{1/u-q},
 &W_q(u,t)&=\Lambda_{u,t}(M_q(u,b)),\label{tail:MW}\\
 D_q(u)&=\prod_{i=q}^{2q-1}(1-iu),
 &V_q(u,t)&=\frac{W_q(u,t)}{c(t)^qD_q(u)}.\label{tail:V}
\end{align}
Empty products equal one.  Expanding in $R-1$ uses at most $q$ factors,
so $M_q$ is a polynomial in $u,b$ despite its displayed form; moreover
$\deg_b M_q\le q$.  Its constant term in $u$ is $c(b)^q$, and hence
\begin{equation}\label{tail:constant-terms}
 W_q(0,t)=c(t)^q,\qquad V_q(0,t)=1.
\end{equation}
Since $q!u^q\binom{1/u-q}{q}=D_q(u)$, the leading comparison in the
physical specialization is
\begin{equation}\label{tail:physical-normalization}
 V_q(1/k,t)=\frac{v_q}{L_q},\qquad
 L_q=2^q\binom{k-q}{q}c(t)^q.
\end{equation}
We use $L_q$ as the comparison sequence.

\subsection{Connected cumulants and logarithmic degree}

For a finite set $S$, let $\Pi(S)$ be its partition lattice and put
$\mu(\pi)=(-1)^{|\pi|-1}(|\pi|-1)!$.  A family of formal units
$(g_I)_{I\subseteq S}$, with $g_\varnothing=1$, has cumulants and
multiplicative interactions
\begin{align*}
 K_J(g)&=\sum_{\pi\in\Pi(J)}\mu(\pi)\prod_{B\in\pi}g_B,\\
 1+\theta_J(g)&=\prod_{I\subseteq J}g_I^{(-1)^{|J|-|I|}}
 \qquad(|J|\ge2).
\end{align*}
We use the following form of the small-cumulant principle
\cite[Proposition 3.3]{DF}.

\begin{lemma}\label{tail:interaction-equivalence}
Suppose $g_I(0)=\prod_{i\in I}g_{\{i\}}(0)$ and all singleton constant
terms are nonzero.  The conditions
\[
 \theta_J(g)=O(u^{|J|-1})\quad(|J|\ge2)
 \quad\hbox{and}\quad
 K_J(g)=O(u^{|J|-1})\quad(|J|\ge2)
\]
are equivalent.
\end{lemma}
\begin{proof}
Dividing $g_I$ by $\prod_{i\in I}g_{\{i\}}$ leaves each interaction
unchanged and multiplies each cumulant by a unit.  Thus assume that all
singletons are identically one.  Multiplicative inversion on subsets
gives
\[
 g_I=\prod_{J\subseteq I,\ |J|\ge2}(1+\theta_J).
\]
On substituting this product into $K_H(g)$, a chosen collection of
hyperedges contributes
$\sum_{\pi\ge\pi_{\rm conn}}\mu(\pi)$, where
$\pi_{\rm conn}$ is its partition into connected components.  This
sum is zero unless the hypergraph is connected, in which case it is
one.  Every connected hypergraph on $H$ satisfies
$\sum_{J\text{ an edge}}(|J|-1)\ge |H|-1$.  The interaction estimates
therefore imply the cumulant estimates.

Conversely, induct on $|H|$.  Keep all proper $g_I$, but replace $g_H$
by $\widetilde g_H=\prod_{J\subsetneq H,\ |J|\ge2}(1+\theta_J)$.
The resulting family has $\theta_H=0$ and satisfies the interaction
estimates by induction.  Its cumulant at $H$ is thus
$O(u^{|H|-1})$.  The difference of the two cumulants is exactly
$g_H-\widetilde g_H$, since every other partition uses proper subsets.
Division by the unit $\widetilde g_H$ proves
$\theta_H(g)=g_H/\widetilde g_H-1=O(u^{|H|-1})$.
\end{proof}

For polynomial inputs define the multilinear cumulant
\begin{equation}\label{tail:joint-definition}
 \kappa_\Lambda(f_1,\ldots,f_r)
 =\sum_{\pi\in\Pi([r])}\mu(\pi)
   \prod_{B\in\pi}\Lambda_{u,t}\left(\prod_{i\in B}f_i\right).
\end{equation}

\begin{lemma}\label{tail:small-Schur}
For $f_1,\ldots,f_r\in\mathbb Q[b]$,
\begin{equation}\label{tail:small-Schur-eq}
 \kappa_\Lambda(f_1,\ldots,f_r)=O(u^{r-1}).
\end{equation}
Every coefficient is a polynomial in $t$, so the assertion includes
$t=0$.  For formal $u$-series $g_i$ with polynomial coefficients it
further implies
\begin{equation}\label{tail:valuation}
 \operatorname{val}_u\kappa_\Lambda(g_1,\ldots,g_r)
 \ge r-1+\sum_i\operatorname{val}_u g_i.
\end{equation}
\end{lemma}
\begin{proof}
First take $f_i=b^{p_i}$ and $t>0$.  Define $\gamma_m(t)$ by
$\log(1-s/t+s^2/t)=\sum_{m\ge1}\gamma_m(t)s^m$.  The product in
\eqref{tail:Lambda} gives the exact formal expansion
\begin{equation}\label{tail:log-weights}
 \log\frac{w_p(u,t)}{t^p}
 =\sum_{m\ge1}\nu^m\left(
 \gamma_m(t)\sum_{i=0}^{p-1}i^m+
 \frac1m\sum_{h=1}^{2p}h^m\right).
\end{equation}
The coefficient of $\nu^m$ is a polynomial in $p$ of degree at most
$m+1$.  Substitution of $\nu=u/(1+u)$ preserves this bound for the
coefficient of $u^m$.

Apply Lemma~\ref{tail:interaction-equivalence} to
$g_I=w_{\sum_{i\in I}p_i}/t^{\sum_{i\in I}p_i}$.  The logarithm of
$1+\theta_J$ is the mixed finite difference of
\eqref{tail:log-weights} with steps $p_i$, $i\in J$.  A mixed difference
of order $|J|$ kills every polynomial of degree less than $|J|$.
Thus $\theta_J=O(u^{|J|-1})$, and the lemma gives
\eqref{tail:small-Schur-eq} for monomials and $t>0$.
The denominator of $w_p$ is independent of $t$, so every coefficient
of the resulting cumulant is polynomial in $t$.  Its vanishing for
all $t>0$ proves the same assertion at $t=0$, also directly visible
from \eqref{tail:t-zero}.  Multilinearity proves the polynomial-input
case and, coefficient by coefficient, \eqref{tail:valuation}.
\end{proof}

To apply this principle to the source, make the formal change of
variable
\begin{equation}\label{tail:change}
 w=zR(z,b),\qquad z=Z(w,b),
\end{equation}
and define
\begin{equation}\label{tail:AD}
 \begin{aligned}
 A(w,b)&=\log R(Z(w,b),b),\\
 \mathcal D(w,b)&=
 \frac{E(Z(w,b),b)}{1+Z(w,b)R_z(Z(w,b),b)/R(Z(w,b),b)}.
 \end{aligned}
\end{equation}
Here $A(0,b)=0$, $[w]A=c(b)$ and $\mathcal D(0,b)=1$.  Formal residue
substitution gives
\begin{equation}\label{tail:residue-change}
 [z^q]ER^{n-q}=[w^q]\mathcal D(w,b)e^{nA(w,b)}.
\end{equation}
Indeed the residue integrand on the left is
$ER^n(zR)^{-q}\,dz/z$, and
$dz/z=(1+zR_z/R)^{-1}\,dw/w$.  Taking $n=1/u$ and replacing $w$ by
$uv$ gives
\begin{equation}\label{tail:conditional-series}
 \sum_{q\ge0}M_q(u,b)\frac{v^q}{q!}
 =\mathcal D(uv,b)e^{A(uv,b)/u}
 =\exp\left(\sum_{r\ge1}s_r(u,b)\frac{v^r}{r!}\right),
\end{equation}
where, if $A=\sum_{r\ge1}\alpha_rw^r$ and
$\log\mathcal D=\sum_{r\ge1}\delta_rw^r$,
\begin{equation}\label{tail:conditional-cumulants}
 s_r(u,b)=r!\bigl(\alpha_r(b)u^{r-1}+\delta_r(b)u^r\bigr),
 \qquad \operatorname{val}_u s_r\ge r-1.
\end{equation}
The substitution makes the exponent independent of $q$.

\begin{theorem}\label{tail:connected-degree}
For each $m\ge1$ there exists a polynomial $\ell_m(q,t)$ in $q$ of
degree at most $m+1$, with coefficients in
$\mathbb Q[t,c(t)^{-1}]$, such that
\begin{equation}\label{tail:degree-eq}
 [u^m]\log V_q(u,t)=\ell_m(q,t)\qquad(q\ge0).
\end{equation}
In particular, these coefficients are continuous for $0\le t\le1/4$,
and
$2\ell_m(q,t)-\ell_m(q-1,t)-\ell_m(q+1,t)$ has degree at most
$m-1$ in $q$.
\end{theorem}
\begin{proof}
For a formal series $K(v)$ with zero constant term, the partition
definition of cumulants gives
\begin{equation}\label{tail:total-identity}
 \log\Lambda_{u,t}(e^{K(v)})
 =\sum_{a\ge1}\frac1{a!}\kappa_\Lambda(K(v),\ldots,K(v)).
\end{equation}
This follows by expanding the exponential and logarithm; each
$v$-coefficient contains only finitely many terms.  Apply it to
$K(v)=\sum_{r\ge1}s_rv^r/r!$.  In a contribution to $[v^n]$, the
input orders $r_1,\ldots,r_a$ sum to $n$.  By
\eqref{tail:valuation} and \eqref{tail:conditional-cumulants}, its
$u$-valuation is at least
$(a-1)+\sum_i(r_i-1)=n-1$.  Thus the ordinary cumulants of the moment
sequence $(W_q)$ satisfy
\begin{equation}\label{tail:small-total}
 \kappa_n^W:=n![v^n]\log\left(\sum_{q\ge0}W_q\frac{v^q}{q!}\right)
 =O(u^{n-1}).
\end{equation}

For every finite set use $g_I=W_{|I|}/c(t)^{|I|}$.  Its constant
terms are one, and \eqref{tail:small-total} gives the cumulant
hypothesis of Lemma~\ref{tail:interaction-equivalence}.  Therefore,
for $r\ge2$,
\[
 \sum_{a=0}^r(-1)^{r-a}\binom ra
 \log\frac{W_a(u,t)}{c(t)^a}=O(u^{r-1}).
\]
At order $u^m$, every forward difference at zero of order $r\ge m+2$
vanishes.  Finite Newton inversion,
$g(q)=\sum_{r=0}^q\binom qr\Delta^rg(0)$, consequently shows that
$[u^m]\log(W_q/c(t)^q)$ agrees, for every nonnegative integer $q$,
with a polynomial of degree at most $m+1$.  Its coefficients belong
to $\mathbb Q[t,c(t)^{-1}]$.  Finally,
\[
 [u^m]\log V_q
 =[u^m]\log(W_q/c(t)^q)+\frac1m\sum_{i=q}^{2q-1}i^m,
\]
and the last power sum has degree at most $m+1$ in $q$.
\end{proof}

\subsection{Uniform factorial estimates and a nonvanishing disk}

For $f\in\mathbb C[b]$ put $\|f\|_1=\sum_p|[b^p]f|$.
Appendix~\ref{tail:appendix} proves the following estimates from a
triangular recurrence for the Schur cumulants.

\begin{proposition}\label{tail:quantitative-inputs}
For $0\le t\le1/4$ the following bounds hold.
\begin{enumerate}
\item If $|\nu|\le1/(8n)$, then for $1\le r\le n$,
\begin{equation}\label{tail:ordinary-bound}
 |\kappa_r^\Lambda(b)|
 \le(r-1)!(64|\nu|)^{r-1}.
\end{equation}
\item The source coefficients in \eqref{tail:conditional-cumulants}
satisfy
\begin{equation}\label{tail:source-bounds}
 \deg\alpha_r,\deg\delta_r\le r,\qquad
 \|\alpha_r\|_1\le\frac{36^r}{r},\quad
 \|\delta_r\|_1\le\frac{56^r}{r},
\end{equation}
and, for $|u|\le1/16$,
\begin{equation}\label{tail:s-bound}
 \|s_r(u,\cdot)\|_1
 \le(r-1)!3136^{r-1}|u|^{r-1}.
\end{equation}
\item If $|u|\le1/(1024n)$ and
$\sum_{i=1}^a\deg f_i\le n$, then
\begin{equation}\label{tail:joint-bound}
 |\kappa_\Lambda(f_1,\ldots,f_a)|
 \le(a-1)!8^n(128|u|)^{a-1}\prod_i\|f_i\|_1.
\end{equation}
If $a\ge2$, a constant input makes this cumulant zero.
\end{enumerate}
Powers with exponent zero are interpreted as one, including at $u=0$.
\end{proposition}

These bounds control the growth with cumulant order in addition to
the valuations of Lemma~\ref{tail:small-Schur}.  They remain uniform
at $t=0$.

\begin{lemma}\label{tail:total-bound}
For $n\ge1$, $0\le t\le1/4$ and $|u|\le1/(1024n)$,
\begin{equation}\label{tail:total-bound-eq}
 |\kappa_n^W|\le8n!26112^{n-1}|u|^{n-1}.
\end{equation}
\end{lemma}
\begin{proof}
Extracting $v^n$ in \eqref{tail:total-identity} gives
\[
 \kappa_n^W=
 \sum_{a=1}^n\frac{n!}{a!}
 \sum_{\substack{r_1+\cdots+r_a=n\\r_i\ge1}}
 \frac{\kappa_\Lambda(s_{r_1},\ldots,s_{r_a})}
 {r_1!\cdots r_a!}.
\]
The polynomial degrees in each summand total at most $n$.  Applying
\eqref{tail:s-bound} and \eqref{tail:joint-bound}, we obtain
\begin{align*}
 |\kappa_n^W|
 &\le n!8^n|u|^{n-1}
 \sum_{a=1}^n\frac{128^{a-1}3136^{n-a}}{a}
 \sum_{\substack{r_1+\cdots+r_a=n\\r_i\ge1}}
 \frac1{r_1\cdots r_a}\\
 &\le n!8^n(128+3136)^{n-1}|u|^{n-1}
 =8n!26112^{n-1}|u|^{n-1}.
\end{align*}
The second inequality bounds the reciprocal factors by one and counts
the compositions by $\binom{n-1}{a-1}$.
\end{proof}

\begin{theorem}\label{tail:disk}
For $q\ge1$ and $0\le t\le1/4$, the rational function $V_q(u,t)$ is
holomorphic and nonzero on a neighborhood of the closed disk
\begin{equation}\label{tail:radius}
 |u|\le\rho_q:=\frac1{2^{25}q}.
\end{equation}
Its logarithm normalized by $\log V_q(0,t)=0$ satisfies
\begin{equation}\label{tail:log-bound}
 |\log V_q(u,t)|<2q\qquad(|u|\le\rho_q).
\end{equation}
\end{theorem}
\begin{proof}
Only moments $w_p$ with $p\le q$ occur, and their denominator factors
are nonzero on this disk.  Thus all functions under consideration are
holomorphic there.  The first conditional cumulant is
$s_1=c(b)(1-u)$, whence
\begin{equation}\label{tail:W1}
 W_1=\frac5{12}(1-u)+\frac t3(1+u)^2.
\end{equation}
For $|u|\le1/16$,
\[
 |W_1-c(t)|\le\frac5{12}|u|+\frac1{12}|u|^2<\frac1{12},
 \qquad |W_1|>\frac13.
\]
Normalize the total cumulants by $\zeta_r=\kappa_r^W/W_1^r$.
Then $\zeta_1=1$, while Lemma~\ref{tail:total-bound} gives
\begin{equation}\label{tail:zeta}
 |\zeta_r|\le(r-1)!(2^{22}|u|)^{r-1}\qquad(2\le r\le q).
\end{equation}
Indeed its immediate bound is
$24r(r-1)!(78336|u|)^{r-1}$, and
$24r\le48^{r-1}$ for $r\ge2$, with
$48\cdot78336<2^{22}$.

Set $Z_m=W_m/W_1^m$ for $0\le m\le q$.  The moment--cumulant
formula, grouping by the block containing one prescribed element,
yields
\begin{equation}\label{tail:partition-recurrence}
 Z_m=Z_{m-1}+\sum_{r=2}^m
 \binom{m-1}{r-1}\zeta_rZ_{m-r},\qquad Z_0=Z_1=1.
\end{equation}
Inductively suppose that the earlier ratios $Z_l/Z_{l-1}$ lie within
$1/2$ of one.  Then
$|Z_{m-r}/Z_{m-1}|\le2^{r-1}$.  Division of
\eqref{tail:partition-recurrence} by $Z_{m-1}$ and
\eqref{tail:zeta} therefore give
\begin{equation}\label{tail:ratio-control}
 \left|\frac{Z_m}{Z_{m-1}}-1\right|
 \le\sum_{h\ge1}(2q\,2^{22}|u|)^h
 \le\sum_{h\ge1}4^{-h}=\frac13<\frac12.
\end{equation}
This proves nonvanishing of every $Z_m$ throughout the closed disk.
Each ratio has the logarithm defined by its series around one, with
modulus at most $\log2<1$.

Using \eqref{tail:V}, we may now write
\[
 \log V_q=q\log\frac{W_1}{c(t)}
 +\sum_{m=1}^q\log\frac{Z_m}{Z_{m-1}}
 -\sum_{i=q}^{2q-1}\log(1-iu).
\]
All terms use the branch vanishing at zero.  Since
$|W_1/c(t)-1|\le2|u|$, the first logarithm has modulus at most
$3|u|$ for $|u|\le1/16$.  Also
$|\log D_q(u)|\le4q^2|u|$ on \eqref{tail:radius}, by summing
$|\log(1-iu)|\le i|u|/(1-i|u|)$.  The displayed logarithm is
therefore bounded by
$q+3q\rho_q+4q^2\rho_q<2q$.  Holomorphy and nonvanishing extend to
a neighborhood of the closed disk because the functions are rational
and the disk is compact.
\end{proof}

Put $C_*=2^{25}$.  Cauchy's estimate gives, for $m\ge1$,
\begin{equation}\label{tail:Cauchy}
 |\ell_m(q,t)|\le2q(C_*q)^m.
\end{equation}
Consequently, if $C_*q/k\le1/2$, then for $M\ge0$,
\begin{equation}\label{tail:remainder}
 \left|\log V_q(1/k,t)-\sum_{m=1}^M
 \frac{\ell_m(q,t)}{k^m}\right|
 \le4q(C_*q/k)^{M+1}.
\end{equation}
For $q=0$ the logarithm and all these coefficients are zero.

\subsection{Discrete curvature}

\begin{lemma}\label{tail:curvature}
Let $q\ge16$ and $1\le m\le\lfloor\sqrt q/4\rfloor$.  Then
\begin{equation}\label{tail:curvature-bound}
 |2\ell_m(q,t)-\ell_m(q-1,t)-\ell_m(q+1,t)|
 \le\frac{432m^2}{q}\left(\frac{3C_*q}{2}\right)^m
 \qquad(0\le t\le1/4).
\end{equation}
\end{lemma}
\begin{proof}
Fix $t$ and put $P(x)=\ell_m(x,t)$, $d=2m$,
$h=\lfloor q/2\rfloor$ and $I=[q-h,q+h]$.
Theorem~\ref{tail:connected-degree} gives $\deg P\le m+1\le d$;
also $h\ge d^2$, $h\ge q/3$ and $h\ge2$.
At every integer $n\in I$, \eqref{tail:Cauchy} gives
\[
 |P(n)|\le3q(3C_*q/2)^m=:G.
\]
The rescaled Markov inequality gives
$\|P'\|_I\le d^2\|P\|_I/h$.  Every point of $I$ is at distance at
most $1/2$ from an integer in $I$, so
$\|P\|_I\le G+d^2\|P\|_I/(2h)$ and hence $\|P\|_I\le2G$.

For $T(\theta)=P(q+h\cos\theta)$, the trigonometric Bernstein
inequality applied twice yields
$\|T'\|_\infty\le d\|P\|_I$ and
$\|T''\|_\infty\le d^2\|P\|_I$.  We use the classical
normalizations in \cite[Theorem 1 and Proposition 2]{Harris}.
When $|x-q|\le1$, we have $|\cos\theta|\le1/2$ and
$\sin\theta\ge\sqrt3/2$ after choosing $0\le\theta\le\pi$.
Differentiating $T$ therefore gives
\[
 |P''(x)|
 \le\frac{\|P\|_I}{h^2}
 \left(\frac43d^2+\frac4{3\sqrt3}d\right)
 \le\frac{2d^2}{h^2}\|P\|_I.
\]
Finally,
\[
 P(q+1)-2P(q)+P(q-1)
 =\int_{-1}^1(1-|y|)P''(q+y)\,dy.
\]
Using $\|P\|_I\le2G$, $h\ge q/3$ and $d=2m$ proves the claimed
constant.
\end{proof}

\begin{proof}[Proof of Theorem~\ref{tail:main}]
Suppose first that $q\ge16$ and put $M=\lfloor\sqrt q/4\rfloor$.
Since $k\ge8192C_*(q+1)$,
$z=3C_*q/(2k)\le3/16384<1/4096$.  Summing
\eqref{tail:curvature-bound} and using
$\sum_{m\ge1}m^2z^m=z(1+z)/(1-z)^3$ shows that the curvature of the
truncated logarithm has absolute value at most
\begin{equation}\label{tail:truncated-curvature}
 \frac{432}{q}\sum_{m\ge1}m^2z^m
 <\frac1{8q}\le\frac1{4(q+1)}.
\end{equation}
For the strict inequality it suffices to substitute $z=1/4096$ and
check $z(1+z)/(1-z)^3<1/3456$.
Using the same $M$ at $q-1,q,q+1$, the total remainder in
\eqref{tail:remainder} is bounded by
\begin{equation}\label{tail:curvature-remainder}
 16(q+1)\left(\frac{C_*(q+1)}k\right)^{M+1}
 \le16(q+1)8192^{-(M+1)}
 \le\frac1{4(q+1)}.
\end{equation}
Indeed, with $n=M+1\ge2$, the definition of $M$ gives
$q+1\le16n^2$.  The inequality $16384n^4\le8192^n$ holds at $n=2$
and then by induction, since $((n+1)/n)^4\le(3/2)^4<8192$.

For $1\le q\le15$, the zeroth-order remainder bound suffices.
Equation~\eqref{tail:remainder} with $M=0$ gives
$|\log V_n(1/k,t)|\le4C_*n^2/k$ for $0\le n\le q+1$.
Consequently the full logarithmic curvature has modulus at most
\begin{equation}\label{tail:small-q}
 \frac{8C_*(2q^2+1)}k
 \le\frac{2q^2+1}{1024(q+1)}
 <\frac1{2(q+1)}.
\end{equation}
Together, \eqref{tail:truncated-curvature}--\eqref{tail:small-q}
bound the curvature of $\log V_q(1/k,t)$ by $1/(2(q+1))$ for every
$q\ge1$ in the theorem.

The baseline in \eqref{tail:physical-normalization} has the exact ratio
\begin{equation}\label{tail:baseline-curvature}
 \frac{L_q^2}{L_{q-1}L_{q+1}}
 =\frac{q+1}{q}\frac{k-q}{k-q+1}
 \frac{(k-2q+2)(k-2q+1)}{(k-2q)(k-2q-1)}
 \ge\frac{q+1}{q}.
\end{equation}
For the last inequality put $s=k-2q>1$ and use
$(s+q)/(s+q+1)\ge s/(s+1)$.  Its logarithm is at least
$1/(q+1)$.  On the positive real specialization the logarithm in
Theorem~\ref{tail:disk} agrees with the real logarithm, because $V_q$
is real, nonvanishing, and equals one at $u=0$.
The normalization relating $v_q$ to $a_{k,j,k-q}$ is independent of
$q$.  Combining \eqref{tail:physical-normalization},
\eqref{tail:baseline-curvature}, and the curvature error proves
\eqref{tail:derivative-margin}.

Apply this estimate for $q=1,\ldots,Q$.  The local transfer
of Theorem~\ref{bridge:local}, using
\[
 g_{k,j,k+m}=\frac{(k+m)!}{m!}\Delta^mF_{k,j}(0),\qquad
 B_{k,j,N}=g_{k,j,N}+g_{k,j,N+1},
\]
then gives \eqref{tail:newton-margin}: the source support makes
$g_{k,j,N}>0$ at all these centers, since $2k-Q\ge k+1$.
The first transferred center is $N=2k-Q$.
\end{proof}

\section{Log-concavity in boundary families}
\label{sec:computational}

We complete Theorem~\ref{intro:families} in two steps.  Near the
balanced index, rational recurrences settle the finite ranges left by
the analytic estimates.  For each fixed Schur index $2\le J\le6$,
we reduce the assertion to polynomial inequalities on a finite
partition of the parameter domain.  The latter reduction, given in
Appendix~\ref{supp:fixed-section}, applies to arbitrarily large degrees.
The corresponding coefficient expansions are collected in Online Resource~2.

\subsection{Finite closure of the central families}

\begin{proposition}[Finite central ranges]\label{comp:central}
With $k=2n+\epsilon$, $j=n-m$ and $\epsilon\in\{0,1\}$, the following
finite assertions hold.
\begin{enumerate}
\item For $m=1$, $1\leq n\leq49$, and for $m=2$,
$2\leq n\leq99$, every original Newton row is positive on its
support from Theorem~\ref{src:gamma-support} and strictly log-concave
at every positive interior center.
\item For $m=3$, $3\leq n\leq124$, and for $m=4$,
$4\leq n\leq179$, every derivative row $(a_{k,j,r})_{r=0}^k$ is
nonnegative and strictly log-concave at $1\leq r\leq k-1$.
Its support is the one in Theorem~\ref{src:shifted-positive}.
\end{enumerate}
\end{proposition}
\begin{proof}
Define
\[
 H_h(d)=[z^{2h}](z/\sin z)^{d+2},\qquad
 S_n(d)=\sum_{h=0}^n\binom nh(2h+1)!d^{2n-2h}H_h(d).
\]
The logarithmic derivative of $z/\sin z$ gives
\begin{equation}\label{comp:H-recurrence}
 H_0=1,\qquad
 H_h=\frac{d+2}{h}\sum_{r=1}^h\tau_rH_{h-r},\qquad
 \tau_r=\frac{(-1)^{r+1}2^{2r-1}B_{2r}}{(2r)!},
\end{equation}
where $B_{2r}$ is an ordinary Bernoulli number; in particular
$\tau_1=1/6$.  The character projection in Section~\ref{sec:near}
gives the following finite combination, up to the indicated positive
scale:
\begin{align}
 P_{2n+\epsilon,n-m}(d)
 &=\lambda_{n,m,\epsilon}\,d^\epsilon
   \sum_{r=0}^m(-1)^{m-r}\binom{m+\epsilon+r}{m-r}
       4^rd^{2r}S_{n-r}(d),\label{comp:raw-source}\\
 \lambda_{n,m,\epsilon}
 &=\frac{\prod_{r=0}^{m-1}(n-r)}
 {4^n n!(n+1)!(n+2)^\epsilon
       \prod_{r=0}^{m-1}(n+\epsilon+2+r)}>0.\notag
\end{align}
The polynomials in \eqref{comp:raw-source} have rational coefficients.
To obtain this expression from the Hahn projection, use
\[
 \binom nh\frac{\fall h\ell}{\fall n\ell}
 =\sum_{r=0}^\ell(-1)^r\binom\ell r\binom{n-r}{h}.
\]
If $u_\ell=(-1)^\ell\fall m\ell
(m+\epsilon+1)^{\overline\ell}/((3/2)^{\overline\ell}\ell!)$,
the Chebyshev expansion used in the character projection gives
\[
 \frac{2m+\epsilon+1}{\epsilon+1}(-1)^r
 \sum_{\ell=r}^m u_\ell\binom\ell r
 =(-1)^{m-r}4^r\binom{m+\epsilon+r}{m-r}.
\]
Substitution proves \eqref{comp:raw-source} with the displayed scale.

For part (i), form $A(d)=\fall{d+1}{k}P(d)$ and compute
\begin{equation}\label{comp:finite-difference}
 B_N=\Delta^NA(0)=\sum_{d=0}^N(-1)^{N-d}\binom NdA(d)
 \qquad(0\leq N\leq2k).
\end{equation}
Since $\deg A\leq2k$, these differences determine the entire row.
For the ranges in part~(i), substitution into
\eqref{comp:H-recurrence}--\eqref{comp:finite-difference} gives
\[
 \begin{aligned}
 B_N&>0 &&\text{on the stated support},\\
 B_N^2-B_{N-1}B_{N+1}&>0
       &&\text{at each positive interior center}.
 \end{aligned}
\]
Each assertion is an inequality between rational numbers; a common
positive denominator reduces it to an integer inequality.

For part~(ii), write the unscaled polynomial on the right of
\eqref{comp:raw-source}, after substituting $d=s+k-1$, as
$\sum_{r=0}^kc_rs^r$.  The required inequalities are
\begin{equation}\label{comp:finite-derivative}
 (r!c_r)^2-(r-1)!(r+1)!c_{r-1}c_{r+1}>0,
 \qquad 1\le r\le k-1.
\end{equation}
The same rational recurrence gives \eqref{comp:finite-derivative}
for every $m,n,\epsilon$ in part~(ii).  The support follows either
from its coefficients or from Theorem~\ref{src:shifted-positive}.
The factor $\lambda_{n,m,\epsilon}$ is positive and therefore does
not affect any of these signs.
The complete primitive integer rows are given in Online Resource~1.
\end{proof}

The analytic ranges begin at $n=50$ for $m=1$ and at
$n=80,125,180$ for $m=2,3,4$, respectively.  Together with
Proposition~\ref{comp:central} and the balanced case proved in
Section~\ref{sec:near}, they cover every imbalance $k-2j\leq9$.
For $m=3,4$, Theorem~\ref{bridge:full} also supplies the factor
\eqref{intro:factor}.

\subsection{A uniform source for the fixed-index kernels}

For the remaining complete families let $J\in\{2,3,4,5,6\}$ be fixed,
and put $m=k-J-1$.  With $S_h(d)=e_h(0,1,\ldots,d)$, extraction of
$[u^J](1-u)C_{k,d}(1,u)$ gives
\begin{equation}\label{comp:fixed-source}
 A_{k,J}(d)=\sum_{s=0}^J(-1)^{J-s}\binom{k-s+1}{J-s}
     d^s\binom{d+1-k+s}{s}S_{k-s}(d).
\end{equation}
Indeed, choose the $k-s$ factors contributing their $(1-u)i$ term in
$i+(d-i)u=(1-u)i+du$, then choose the remaining $s$ positions.
For $\mathcal N(\sum b_N\binom dN)=\sum b_Nt^N$, one has
\[
 \mathcal N(df)=M\mathcal Nf,\quad M=t(1+t)\partial_t+t,
 \qquad\mathcal NS_{h+1}=L\mathcal NS_h,\quad L=t(M+1).
\]
The second identity follows from
$\Delta S_{h+1}(d)=(d+1)S_h(d)$ and $S_{h+1}(0)=0$.

Starting with $\mathcal NS_1=t(1+t)$, write
\[
 \begin{gathered}
 \mathcal NS_{m+1}=t^{m+1}(1+t)R_m(t),\qquad R_0=1,\\
 R_{m+1}=(m+2+(m+3)t)R_m+t(1+t)R'_m.
 \end{gathered}
\]
Induction at consecutive roots shows that $R_m$ has $m$ simple roots
in $(-1,0)$: the recurrence alternates in sign at its old roots, and
$R_m(0)=(m+1)!$, $R_m(-1)=(-1)^m$ provide the two endpoint signs.
Writing
\[
 R_m(t)=\sum_{a=0}^m w_at^a(1+t)^{m-a}
\]
therefore gives a polynomial
$\sum_aw_az^a=(1-z)^mR_m(z/(1-z))$ with negative real roots and
positive coefficients.  Its coefficients $w_a$ are log-concave by
Newton's inequalities.

For independent $a,b$, let $k=a+b+J+1$ and set
\[
 \widetilde M=t(1+t)\partial_t+2a+b+1+(2a+2b+3)t,
 \qquad\widetilde L=t(\widetilde M+1).
\]
Define the degree-$2J$ polynomial
\begin{equation}\label{comp:kernel-definition}
\begin{split}
 P^{(J)}_{a,b}(t)=J!\sum_{s=0}^J
 &\frac{(-1)^{J-s}}{s!}\binom{k-s+1}{J-s}\widetilde M^s\\
 &\cdot\prod_{v=0}^{s-1}(\widetilde M-(k-s-1+v))
           \widetilde L^{J-s}1.
\end{split}
\end{equation}
Operators act from right to left.  Conjugating $M,L$ by
$t^{m+1+a}(1+t)^{b+1}$ in \eqref{comp:fixed-source} gives
\begin{equation}\label{comp:kernel-output}
 \sum_NB_{k,J,N}t^N=t^{k-1}\sum_{a=0}^mw_aF_a(t),\qquad
 F_a(t)=\frac{t^{a-J+1}(1+t)^{m-a+1}}{J!}P^{(J)}_{a,m-a}(t).
\end{equation}
Thus the kernel is specified directly by the original Chern product
and the two Newton operators.

\subsection{Centered prefixes}

Put $K_{r,a}=[t^r]F_a(t)$.  The sufficient LC-positivity criterion
of Wang--Yeh \cite{WY} asks for
\begin{equation}\label{comp:prefix}
 \mathcal P_{r,S,a}:=\sum_{c=a}^{S-a}
 \bigl(K_{r,c}K_{r,S-c}-K_{r-1,c}K_{r+1,S-c}\bigr)\geq0
\end{equation}
on every centered interval $0\leq a\leq S-a\leq m$.
For completeness, pair the two orders off the diagonal and denote the
paired summands by $T_a$, counting a central diagonal once.  If
$z_a=w_aw_{S-a}$, then $z_a$ is nondecreasing from the outside toward
the center.  With $C_a=\sum_{c=a}^{\lfloor S/2\rfloor}T_c$ and
$l=\max(0,S-m)$, summation by parts gives
\[
 \sum_{a=l}^{\lfloor S/2\rfloor}T_az_a
 =C_lz_l+\sum_{a=l+1}^{\lfloor S/2\rfloor}C_a(z_a-z_{a-1})\geq0.
\]
Summing over $S$ proves the output Tur\'an inequality.

\begin{proposition}[Fixed Schur indices]\label{comp:fixed}
For each $J\in\{2,3,4,5,6\}$, the original row
$(B_{k,J,N})_N$ is log-concave for every $k\geq2J$.
\end{proposition}
\begin{proof}
Appendix~\ref{supp:fixed-section} establishes the requisite
polynomial inequalities on the ranges in Table~\ref{tab:fixed-ranges}.
An interval whose
upper endpoint is smaller than its lower endpoint is empty.

\begin{table}[htbp]
\caption{Ranges covered by the fixed-index kernel argument.}
\label{tab:fixed-ranges}
\centering\small\setlength{\tabcolsep}{4pt}
\begin{tabular}{@{}ccccc@{}}
\toprule
$J$ & Kernel range & Interior centers $r$ & Upper centers $r$ & Extra degrees\\
\midrule
2 & $m\geq1$ & $1\leq r\leq m-1$ & $m\leq r\leq m+3$ & none\\
3 & $m\geq2$ & $1\leq r\leq m-2$ & $m-1\leq r\leq m+4$ & none\\
4 & $m\geq3$ & $1\leq r\leq m-3$ & $m-2\leq r\leq m+5$ & none\\
5 & $m\geq5$ & $1\leq r\leq m-4$ & $m-3\leq r\leq m+6$ & $k=10$\\
6 & $m\geq7$ & $1\leq r\leq m-6$ & $m-5\leq r\leq m+7$ & $k=12,13$\\
\bottomrule
\end{tabular}
\end{table}

First, the kernel coefficients are nonnegative and log-concave on
the stated domains, and every apparent negative power in
\eqref{comp:kernel-output} cancels.  Each column is consequently
log-concave.  For the interior centers, common-binomial identities
express \eqref{comp:prefix} as the sum of a telescoping boundary term
and a remainder.  The affine substitutions in
Appendix~\ref{supp:domains} give nonnegative coefficient expansions
of both terms on every nonexceptional domain.  The intervals $[0,1]$
for $J=4$ and $[0,1],[0,2]$ for $J=5,6$ require the complete prefix,
rather than separate signs for its two terms.  Their full sums also
have nonnegative expansions.  For the upper centers, the complete
ordered sum in \eqref{comp:prefix} has a nonnegative expansion in
ordinary powers and binomial polynomials of the free gap parameter.
All of these bases are nonnegative on their stated integer domains.

All nonzero columns have a common right endpoint $m+J+2=k+1$.
The actual weights are positive, and one of the initial columns is
positive over the whole required interval.  At the only exceptional
column boundary $J=4,m=3,a=0$, column $a=1$ supplies that interval.
Thus the output has no internal zeros.  Table~\ref{tab:fixed-ranges} covers
every center $1\leq r\leq m+J+1$; the endpoint inequalities follow
from zero extension.  The preceding Abel argument proves
log-concavity for every degree in the kernel range.

The three remaining pairs $(J,k)=(5,10),(6,12),(6,13)$ are all
balanced: $J=\lfloor k/2\rfloor$.  Their assertions follow from
Theorem~\ref{near:balanced} and Theorem~\ref{bridge:full}.  This
covers every $k\geq2J$.
\end{proof}

Propositions~\ref{comp:central} and~\ref{comp:fixed}, the balanced
and first-mode arguments in Section~\ref{sec:near}, and the known
$j=0,1$ cases prove Theorem~\ref{intro:families}.  The positivity
and support of every row also follow independently from
Theorem~\ref{src:gamma-support}.

\section{The remaining Newton problem}
\label{sec:discussion}

The inverse-series formula establishes nonnegativity and exact support
for all Newton rows.  Log-concavity is proved in the uniform regions of
Theorem~\ref{intro:regions} and the boundary families of
Theorem~\ref{intro:families}.  For the unrestricted problem, the transfer
in Section~\ref{sec:bridge} suggests the following stronger statement.

\begin{conjecture}\label{discussion:derivative}
For every $k\ge1$ and $0\le j\le\lfloor k/2\rfloor$, the sequence
\[
 a_{k,j,r}=[s_{(k-j,j)}][z^k]e^{\Ccal(z)}H(z)^r,
 \qquad 0\le r\le k,
\]
is log-concave in $r$.
\end{conjecture}

Here $[s_\lambda]$ denotes Schur coefficient extraction.
Theorem~\ref{src:shifted-positive} gives positivity and interval support,
so the conjecture concerns only the adjacent inequalities.
By Theorem~\ref{bridge:full}, it would imply the original Newton
conjecture, together with the stronger factor \eqref{intro:factor}
throughout the positive interior.

The real-rootedness established in Section~\ref{bern:section} controls
the binary variables at a fixed shift.  It therefore proves
ultra-log-concavity in the Schur index, while the conjecture above
requires an inequality in the exponent of $H$.  Positivity of the
quotient alone also does not suffice.  For example, at $k=2$, the
polynomial $F(s)=s^2+s+100$ gives
\[
 A(d)=d(d+1)F(d-1),\qquad (B_1,B_2,B_3,B_4)=(200,212,36,24),
\]
and $36^2<212\cdot24$.  Although this quotient does not arise from
the Chern construction, it shows why a positivity argument requires
an additional inequality before the Newton transfer can be applied.

Two quantitative questions remain within the methods developed here.
Sharper estimates for the Hahn projection may enlarge the range in
which the shifted zeros lie in the required sector.  Stronger cumulant
bounds may enlarge the terminal region, since the reduction of
logarithmic degree already holds at every order.  The fixed-index
argument suggests a separate algebraic question: whether the
centered-prefix inequalities admit a proof uniform in the Schur index.
Such a proof would remove the need to treat exceptional small prefixes
individually.

\begin{appendices}
\renewcommand{\theHtable}{\thesection.\arabic{table}}
\section{The complex saddle estimate for the linear range}
\label{near:linear-proof}

We prove Theorem~\ref{near:linear} by a saddle estimate uniform in
complex $d$ outside the stated zero disk.  Introducing an auxiliary
complex frequency gives control as $\omega/n$ tends to zero.

\subsection{An integral representation}

Let
\[
 \Phi(w)=\log\frac{w}{\sin w},\qquad \Phi(0)=0,\qquad |w|<\pi,
\]
with the analytic branch at zero.  The character projection of
\eqref{near:source} is
\[
 P_{k,j}(d)=\frac2\pi[t^k]\int_0^\pi
 \sin\theta\sin(\omega\theta)
 e^{dt\cos\theta+(d+2)\Phi(t\sin\theta)}\,d\theta.
\]
The parity of $k$ and $\omega=2m+\epsilon+1$ allows reflection in
$\pi/2$ to halve the interval.  A subsequent reflection in zero
cancels the odd cosine-character term, giving
\begin{equation}\label{near:oscillatory-source}
 P_{k,j}(d)=\frac{2}{\pi i}[t^k]
 \int_{-\pi/2}^{\pi/2}\sin\theta\,e^{i\omega\theta}
 e^{dt\cos\theta+(d+2)\Phi(t\sin\theta)}\,d\theta.
\end{equation}
To exclude zeros with $|d|\ge7n/5$, put
\[
 \begin{gathered}
 z=\frac nd,\quad 0<|z|\le\frac57,\qquad
 \varkappa=\frac kn,\quad\eta=\frac2n,\quad\rho=\frac\omega n,\\
 \gamma=1+\eta z,\qquad c=\frac\gamma z,\qquad t_* =\varkappa z.
 \end{gathered}
\]
Define
\begin{equation}\label{near:phase}
 \Psi(\varphi,\theta)
 =\varkappa e^{i\varphi}\cos\theta
  +c\Phi(\varkappa z e^{i\varphi}\sin\theta)
  -i\varkappa\varphi+i\rho\theta.
\end{equation}
Since $\Phi(w)=O(w^2)$, the term involving $c$ extends analytically
to $z=0$.  Applying Cauchy's formula on $t=t_*e^{i\varphi}$ to
\eqref{near:oscillatory-source} gives
\begin{equation}\label{near:double-integral}
 P_{k,j}(d)=\frac{t_*^{-k}}{\pi^2 i}I(\rho),\qquad
 I(\rho)=\int_{-\pi}^{\pi}\int_{-\pi/2}^{\pi/2}
       \sin\theta\,e^{n\Psi(\varphi,\theta)}\,d\theta\,d\varphi.
\end{equation}

Fix $n,\epsilon,z$, regard $\rho$ as an independent complex
parameter, and put
\[
                 \delta=10^{-12},\qquad \sigma=|\rho|\le\delta.
\]
The integral is defined for every complex $\rho$.  We recover the
Schur quotient from \eqref{near:double-integral} at the integer
frequency $n\rho=\omega$ after establishing the uniform estimate.

\subsection{Uniform bounds for the phase}

Euler's product gives
\[
 \Phi(w)=\sum_{r\ge1}\frac{\tau_r}{r}w^{2r},\qquad
 \tau_1=\frac16,\quad\tau_2=\frac1{90},\quad
 0<\tau_r\le\frac1{90\,9^{r-2}}\quad(r\ge2).
\]
Using $(\pi^2\nu^2)^{-1}<1/9$ in the last estimate gives, for
$|w|\le3/2$,
\begin{equation}\label{near:Phi-global}
             |\Phi(w)|\le\frac{|w|^2}{5},\qquad
             |\Phi'(w)|\le\frac25|w|.
\end{equation}
The derivative is bounded by
$|w|/3+|w|^3/[45(1-|w|^2/9)]\le2|w|/5$; integration along a
radius gives the first inequality.  On $|w|\le1/2$,
\begin{equation}\label{near:Phi-local}
 |\Phi'(w)-w/3|\le\frac{|w|^3}{40},\qquad
 |\Phi''(w)-1/3|\le\frac{|w|^2}{10},\qquad
 |\Phi''(w)|\le\frac25.
\end{equation}
Differentiating the geometric majorant gives the middle bound:
\[
 |\Phi''(w)-1/3|
 \le\frac{|w|^2}{45}
       \frac{3-|w|^2/9}{(1-|w|^2/9)^2}
 <\frac{|w|^2}{10};
\]
the first bound follows from the corresponding majorant for $\Phi'$.

For $n\ge100$, the following bounds hold uniformly on the prescribed
$z$-disk:
\begin{equation}\label{near:saddle-parameters}
 \begin{gathered}
 2\le\varkappa\le\frac{201}{100},\qquad
 |\varkappa z|\le\frac{201}{140}<\frac32,\\
 |\gamma|\le\frac{51}{50},\qquad
 M:=\varkappa^2|z+\eta z^2|<\frac{147}{50}.
 \end{gathered}
\end{equation}
The last inequality uses $(201/100)^2(5/7+1/98)<147/50$.
We also require the bound
\begin{equation}\label{near:third-derivatives}
 |\partial_\varphi^a\partial_\theta^b\Psi|<600
 \qquad(a+b=3),\qquad |\varphi|,|\theta|\le1/10.
\end{equation}
The domains here and below are complex polydisks.  On
$|\varphi|,|\theta|\le1/2$, the inequalities $e^{1/2}<5/3$,
$\cosh(1/2)<8/7$, and $\sinh(1/2)<8/15$ bound the modulus of
the argument of $\Phi$ by $3/2$.  Thus
\eqref{near:Phi-global}--\eqref{near:saddle-parameters} give
\[
 |\Psi|\le
 \frac{201}{100}\frac53\frac87
 +\frac{147}{250}\left(\frac53\right)^2
                       \left(\frac8{15}\right)^2
 +\frac{201}{200}+\frac\sigma2<6.
\]
A polydisk of radius $2/5$ about any point of the smaller domain
lies in the larger domain.  Cauchy's estimate therefore bounds each
derivative of total order three by $3!\,6/(2/5)^3=1125/2<600$.

\subsection{The critical point}

Set
\[
 u=e^{i\varphi}\cos\theta,\qquad
 v=e^{i\varphi}\sin\theta.
\]
Differentiation of \eqref{near:phase} gives the critical-point equations
\begin{equation}\label{near:saddle-equations}
 u=1-\gamma v\Phi'(\varkappa zv),\qquad
 v-\gamma u\Phi'(\varkappa zv)=\frac{i\rho}{\varkappa}.
\end{equation}
Eliminating $u$ gives the scalar odd equation
\begin{align}
 Dv+\mathcal R(v)&=\frac{i\rho}{\varkappa},\qquad
 D=1-\frac{\varkappa z\gamma}{3},\label{near:saddle-scalar}\\
 \mathcal R(v)
 &=-\gamma\{\Phi'(\varkappa zv)-\varkappa zv/3\}
       +\gamma^2v\Phi'(\varkappa zv)^2.\nonumber
\end{align}
By \eqref{near:saddle-parameters}, $|D-1|<1/2$, so
$\Re D>1/2$ and $|D|<3/2$.  For $|v|\le1/100$,
\eqref{near:Phi-global} and \eqref{near:Phi-local} yield
\begin{equation}\label{near:R-bounds}
                    |\mathcal R(v)|\le|v|^3,
          \qquad |\mathcal R'(v)|\le2|v|^2.
\end{equation}
The coefficient bounds are
\[
 \frac{51}{50}\frac1{40}\left(\frac32\right)^3
 +\left(\frac{51}{50}\right)^2
        \left(\frac25\frac32\right)^2<1
\]
and
\[
 \frac{51}{50}\frac1{10}\left(\frac32\right)^3
 +\left(\frac{51}{50}\right)^2
  \left\{\left(\frac25\frac32\right)^2
        +2\left(\frac25\frac32\right)\frac32\frac25\right\}<2,
\]
respectively.

Let $v_0=i\rho/(\varkappa D)$.  For $\rho\ne0$,
$\sigma/4<|v_0|<\sigma$.  On $|v|\le2\sigma$, the map
$v\mapsto v_0-\mathcal R(v)/D$ maps the disk into itself and
has Lipschitz constant at most $16\sigma^2<1/2$.
Its unique fixed point satisfies
\begin{equation}\label{near:v-bounds}
 |v-v_0|\le16\sigma^3,\qquad
 |v|\le2\sigma,\qquad |v|\ge\sigma/5.
\end{equation}
At $\rho=0$, the fixed point is $v=0$.  On the fixed disk
$|v|\le2\delta$, the iterates converge uniformly in $\rho$.
The solution is therefore holomorphic on a neighborhood of
$|\rho|\le\delta$.

Define $u$ by the first equation of \eqref{near:saddle-equations}.
Then
\[
 |u-1|\le3\sigma^2,\qquad |u^2+v^2-1|\le11\sigma^2.
\]
The branches near $1$ and $0$ in
\begin{equation}\label{near:saddle-location}
 \varphi_s=\frac{\log(u^2+v^2)}{2i},\qquad
 \theta_s=\arctan(v/u)
\end{equation}
are well defined and recover $u,v$, hence a critical point of the
phase.  The logarithm and arctangent series give
\begin{equation}\label{near:saddle-location-bounds}
 |\varphi_s|\le10\sigma^2,\qquad |\theta_s|\le4\sigma,\qquad
 \frac\sigma{10}\le|a_s|\le3\sigma,\qquad |b_s|\le2,
\end{equation}
where $a_s=\sin\theta_s$ and $b_s=\cos\theta_s$.
The amplitude bounds also follow from
$a_s=v/\sqrt{u^2+v^2}$ and $b_s=u/\sqrt{u^2+v^2}$.
All these quantities depend holomorphically on $\rho$.

Let $K=-\Psi''(\varphi_s,\theta_s)$, a complex symmetric matrix.
At the origin the comparison matrix is
\[
 K_0=\begin{pmatrix}
 \varkappa&0\\
 0&\varkappa-\varkappa^2(z+\eta z^2)/3
 \end{pmatrix}.
\]
The real part of the comparison matrix is at least the identity, and its operator
norm is less than $3$.  Equations~\eqref{near:third-derivatives}
and \eqref{near:saddle-location-bounds} bound each entry of
$K-K_0$ by $3000\sigma$, giving
\begin{equation}\label{near:Hessian}
 \Re(X^TKX)\ge\frac9{10}|X|^2\quad(X\in\mathbb R^2),
                  \qquad \|K\|<4.
\end{equation}

\subsection{Contour deformations and the Gaussian comparison}

Translate the $\varphi$ contour to $[-\pi,\pi]+\varphi_s$.
Since $n\varkappa=k$ is an integer, the integrand in
\eqref{near:double-integral} is $2\pi$-periodic in $\varphi$,
and the connecting sides cancel for every complex $n\rho$.
Then translate the $\theta$ segment to
$[-\pi/2,\pi/2]+\theta_s$ and retain its two connecting segments.
Along both deformations,
\begin{equation}\label{near:no-singularities}
 |\varkappa z e^{i\varphi}\sin\theta|
 \le\frac{201}{140}\left(\frac{101}{100}\right)^2
       <\frac32<\pi,
\end{equation}
since $e^{|\Im\varphi_s|}$ and $\cosh|\Im\theta_s|$ are each
less than $101/100$.  The integrand remains analytic throughout
the deformation, with the same branch of $\Phi$.

Put $\Psi_s=\Psi(\varphi_s,\theta_s)$ and use real coordinates
$X=(x,y)$ on the translated rectangle.  On the square
\[
 \mathcal C=\{(x,y):|x|,|y|\le r_0\},\qquad r_0=10^{-4},
\]
Taylor expansion at the critical point gives
\begin{equation}\label{near:Taylor}
 \Psi(\varphi_s+x,\theta_s+y)
       =\Psi_s-\tfrac12X^TKX+E(X),
\end{equation}
where
\begin{equation}\label{near:Taylor-error}
 |E(X)|\le100(|x|+|y|)^3\le300|X|^3,\qquad
 |E(X)|\le\frac1{25}|X|^2.
\end{equation}
The first bound follows from \eqref{near:third-derivatives}, including
the Taylor factor $1/3!$.  The second follows from
$(|x|+|y|)^3\le(2r_0)\,2|X|^2$ and $400r_0=1/25$.
The critical-point equations eliminate the linear term in
\eqref{near:Taylor}.

For $\rho\ne0$, use the Gaussian comparison integral
\begin{equation}\label{near:Gaussian}
 I_G=a_s e^{n\Psi_s}\int_{\mathbb R^2}
                   e^{-nX^TKX/2}\,dX
       =a_s e^{n\Psi_s}\frac{2\pi}{n\sqrt{\det K}}.
\end{equation}
The formula for real positive definite matrices extends
holomorphically to the convex domain of complex symmetric matrices
with positive real part.  The square root is chosen positive on the
real positive definite matrices.  By \eqref{near:Hessian},
$|\sqrt{\det K}|\le\|K\|<4$, and hence
\begin{equation}\label{near:Gaussian-lower}
 |I_G|\ge\frac{|a_s|}{n}e^{n\Re\Psi_s}
          \ge\frac\sigma{10n}e^{n\Re\Psi_s}.
\end{equation}

\subsection{The central error on the parameter circle}

In this and the next subsection, assume $|\rho|=\delta$, so
$\sigma=\delta>0$.  To keep track of the amplitude $a_s$, use the
sine addition formula for real $y$:
\begin{equation}\label{near:sine-amplitude}
 \sin(\theta_s+y)=a_s+b_sy+e_a(y),\qquad
 |e_a(y)|\le\tfrac12|a_s|y^2+\tfrac13|y|^3.
\end{equation}
Central symmetry under $X\mapsto-X$ makes the integral of $b_sy$
against the centered Gaussian vanish on both $\mathcal C$ and
$\mathbb R^2$.  The relative estimate near zero frequency relies on
this cancellation.

Remove the common factor $e^{n\Psi_s}$.  The inequality
$|e^{nE}-1|\le n|E|e^{n|E|}$, together with
\eqref{near:Hessian}--\eqref{near:sine-amplitude}, bounds the error
over $\mathcal C$ by
\begin{equation}\label{near:central-integral-error}
 \int_{\mathbb R^2}e^{-(2/5)n|X|^2}
 \left\{300n|a_s||X|^3+600n|X|^4
              +\tfrac12|a_s||X|^2+\tfrac13|X|^3\right\}\,dX.
\end{equation}
The exponential decay has coefficient at least
$9/20-1/25=41/100>2/5$.
Polar coordinates, or
\[
 \int_{\mathbb R^2}|X|^p e^{-an|X|^2}\,dX
          =\frac{\pi\,\Gamma((p+2)/2)}{(an)^{(p+2)/2}},
\]
give
\[
 \int |X|^2e^{-(2/5)n|X|^2}dX<20n^{-2},\qquad
 \int |X|^3e^{-(2/5)n|X|^2}dX<50n^{-5/2},
\]
and
$\int |X|^4e^{-(2/5)n|X|^2}dX<100n^{-3}$.
Dividing \eqref{near:central-integral-error} by $|a_s|/n$ and
using $|a_s|\ge\delta/10$ bounds the relative error by
\begin{equation}\label{near:central-relative}
 \frac{15000}{\sqrt n}+\frac{600000}{n\delta}
       +\frac{10}{n}+\frac{200}{\delta n^{3/2}}<10^{-8}.
\end{equation}
Each term decreases with $n$, and the bound holds at $n=10^{40}$.

Extending the constant-amplitude Gaussian from $\mathcal C$ to
$\mathbb R^2$ has relative error at most
\begin{equation}\label{near:Gaussian-tail}
 \frac{20\pi}{9}\exp\{-9nr_0^2/20\}
       <10\exp\{-n/(3\cdot10^8)\}<10^{-8}.
\end{equation}
Indeed, the complement of $\mathcal C$ lies in $|X|\ge r_0$,
where the Gaussian integral is bounded by
$(20\pi/(9n))e^{-9nr_0^2/20}$.  Division by
\eqref{near:Gaussian-lower} gives the first inequality.  Since the
linear-amplitude integral vanishes on both domains, the error
retains the factor $|a_s|$.

\subsection{Global tails and connecting segments}

On the original real rectangle write $c_0=\cos\theta\in[0,1]$.
Equations~\eqref{near:Phi-global} and
\eqref{near:saddle-parameters} imply
\[
 \Re\Psi(\varphi,\theta)
       \le\varkappa\cos\varphi\,c_0
              +\frac35(1-c_0^2)+2\sigma.
\]
The frequency term satisfies
$\Re(i\rho\theta)\le\sigma\pi/2<2\sigma$ for complex $\rho$.
After allowing for this term, the gap from
$\varkappa=\Psi(0,0)$ is at least
\begin{equation}\label{near:real-gap}
             \frac45(1-c_0)+2c_0(1-\cos\varphi).
\end{equation}
For $|t|\le\pi$ one has $1-\cos t\ge t^2/5$; for example,
$\sin u\ge2u/\pi$ on $[0,\pi/2]$ gives the stronger constant
$2/\pi^2>1/5$.  If $|\theta|\ge r_0$, the first term of
\eqref{near:real-gap} is at least $4r_0^2/25$.
Otherwise $c_0>1/2$, and $|\varphi|\ge r_0$ makes the second
term at least $r_0^2/5$.  Outside $\mathcal C$, the gap on the
real rectangle is therefore at least
$r_0^2/10-2\sigma=10^{-9}-2\sigma$.

On the thin complex strips used for the deformations, direct
differentiation of \eqref{near:phase} gives
\begin{equation}\label{near:strip-derivatives}
                     |\Psi_\varphi|<6,\qquad |\Psi_\theta|<4.
\end{equation}
The three contributions to the first bound sum to at most
\[
 \frac{201}{100}\left(\frac{101}{100}\right)^2
 +\frac25\frac{147}{50}\left(\frac{101}{100}\right)^4
 +\frac{201}{100}<6.
\]
For the second bound, replacing the last summand by $\sigma$
gives a quantity less than $4$.  Translation therefore changes
$\Psi$ by at most $6|\varphi_s|+4|\theta_s|\le17\sigma$;
in particular, $|\Psi_s-\varkappa|\le17\sigma$.
On the translated rectangle outside $\mathcal C$, this yields
\begin{equation}\label{near:translated-gap}
 \Re\Psi(\varphi_s+x,\theta_s+y)
 \le\Re\Psi_s-10^{-9}+36\sigma
 \le\Re\Psi_s-\frac1{2\cdot10^9}.
\end{equation}
The rectangle has area $2\pi^2<20$, and $|\sin\theta|<2$.
Division by \eqref{near:Gaussian-lower} bounds this global tail
relative to $I_G$ by
\begin{equation}\label{near:global-relative}
              \frac{400n}{\delta}
              e^{-n/(2\cdot10^9)}<10^{-8}.
\end{equation}

The two connecting $\theta$ segments have the form
$\theta=\pm\pi/2+u\theta_s$, $0\le u\le1$, with
$\varphi=x+\varphi_s$.  There,
$|\cos\theta|\le5\sigma$ and
$|\sin\theta|<101/100$.  The four terms in
\eqref{near:phase} have respective real-part upper bounds
\[
 11\sigma,\qquad
 \frac{147}{250}\left(\frac{101}{100}\right)^4,\qquad
 21\sigma^2,\qquad 2\sigma,
\]
so that $\Re\Psi<7/10$.
Each segment has length at most $4\sigma$, whereas
$\Re\Psi_s\ge2-17\sigma>19/10$.
Their combined contribution, relative to
\eqref{near:Gaussian-lower}, is at most
\begin{equation}\label{near:connecting-relative}
                           1200n e^{-n}<10^{-8}.
\end{equation}
The exponential comparisons in \eqref{near:Gaussian-tail},
\eqref{near:global-relative}, and \eqref{near:connecting-relative}
follow from $e^x\ge x^4/24$: the resulting rational-power bounds
decrease with $n$ and are smaller than $10^{-8}$ at $n=10^{40}$.

\subsection{Continuation to zero frequency}

On the parameter circle $|\rho|=\delta$, the four relative errors
in \eqref{near:central-relative}, \eqref{near:Gaussian-tail},
\eqref{near:global-relative}, and \eqref{near:connecting-relative}
sum to less than $4\cdot10^{-8}<1$.  Holomorphic continuation of
the relative error extends this estimate through $\rho=0$.

The integral $I(\rho)$ in \eqref{near:double-integral} is entire
in $\rho$.  Reflection $\theta\mapsto-\theta$ on the original
rectangle shows that it is odd.
Uniqueness in \eqref{near:saddle-scalar} implies that $v(\rho)$
is odd and $u(\rho)$ is even.  The branches in
\eqref{near:saddle-location} therefore make $\varphi_s$ even
and $\theta_s$ odd.  It follows that $\Psi_s$ and $\det K$
are even, while $a_s=\sin\theta_s$ is odd.  Moreover,
\[
          \frac{a_s(\rho)}\rho\longrightarrow
                   \frac{i}{\varkappa D}\ne0
                   \qquad(\rho\longrightarrow0).
\]
By \eqref{near:saddle-location-bounds}, zero is the only zero of
$a_s$ in the disk.  Since $\det K$ is holomorphic and nonvanishing,
its square root has the Gaussian branch of \eqref{near:Gaussian}.
The value at zero and uniqueness of continuation make this branch
even.  Consequently
\begin{equation}\label{near:removable-error}
 \mathcal E(\rho)=
 \frac{I(\rho)e^{-n\Psi_s(\rho)}n\sqrt{\det K(\rho)}}
      {2\pi a_s(\rho)}-1
\end{equation}
is holomorphic on a neighborhood of $|\rho|\le\delta$,
including its removable singularity at zero.
The maximum modulus principle and the circle estimate give
\[
               |\mathcal E(\rho)|<4\cdot10^{-8}<1
                      \qquad(|\rho|\le\delta).
\]
For every $0<|\rho|\le\delta$ the Gaussian factor is nonzero,
so $I(\rho)\ne0$.

Set $\rho=\omega/n$.  The estimates hold uniformly for complex
$z$ with $0<|z|\le5/7$, so \eqref{near:double-integral} shows that
$P_{k,j}(d)$ has no zero with $|d|\ge7n/5$.  This proves
Theorem~\ref{near:linear}.
For its simpler sufficient range, $2m+\epsilon+1\le2m+2$
is at most $10^{-12}n$ whenever
$m\le n/(4\cdot10^{12})$ and $n\ge10^{40}$.

\section{Uniform cumulant estimates}
\label{tail:appendix}

We prove Proposition~\ref{tail:quantitative-inputs} by working in
finite coefficient ranges.  The moment series is used formally
throughout.

\subsection{A triangular recurrence for Schur cumulants}

Regard $\nu$ as the parameter and write $w_p=w_p(\nu,t)$ for the
product in \eqref{tail:Lambda}.  Set
\[
 K(z)=\sum_{p\ge0}w_p\frac{z^p}{p!},\qquad
 \kappa_r^\Lambda(b)=r![z^r]\log K(z).
\]
The weights satisfy the first-order recurrence
\begin{equation}\label{tail:weight-recurrence}
 (1-(2p+1)\nu)(1-(2p+2)\nu)w_{p+1}
 =(t-p\nu+p^2\nu^2)w_p.
\end{equation}
For $\vartheta=z\partial_z$, this is equivalent coefficientwise to
\[
 (1-\nu-2\nu\vartheta)(1-2\nu-2\nu\vartheta)K'
 =(t-\nu\vartheta+\nu^2\vartheta^2)K,
\]
or, after expansion,
\begin{equation}\label{tail:moment-ODE}
\begin{split}
 (1-3\nu+2\nu^2)K'&+(-4\nu+10\nu^2)zK''
  +4\nu^2z^2K'''\\
 &=tK+(-\nu+\nu^2)zK'+\nu^2z^2K''.
\end{split}
\end{equation}
Divide by $K$ and substitute
\[
 \partial_z\log K(z)=\sum_{h\ge0}C_h(\nu,t)(\nu z)^h.
\]
Coefficient extraction gives
$C_0=t/((1-\nu)(1-2\nu))$ and, for $h\ge1$,
\begin{equation}\label{tail:triangular}
\begin{split}
 &(1-(2h+1)\nu)(1-(2h+2)\nu)C_h\\
 &\quad=(-1+h\nu)C_{h-1}
  +[4-(6h+4)\nu]\sum_{i+j=h-1}C_iC_j\\
 &\qquad\quad+\sum_{i+j=h-2}C_iC_j
  -4\sum_{i+j+l=h-2}C_iC_jC_l.
\end{split}
\end{equation}
Empty sums are zero.  Performing the substitution first for
$\nu\ne0$, the recurrence extends each $C_h$ regularly to
$\nu=0$ and gives
\begin{equation}\label{tail:kappa-C}
 \kappa_{h+1}^\Lambda(b)=h!\nu^hC_h.
\end{equation}
The divisibility by $\nu^h$ can also be read inductively from the
recurrence for $[z^h]\partial_z\log K$, before division by $\nu^h$.

\begin{proof}[Proof of Proposition~\ref{tail:quantitative-inputs}(i)]
The estimate holds for complex $t$ with $|t|\le1/4$.
If $0\le h\le n-1$ and $|\nu|\le1/(8n)$, the two factors on
the left of \eqref{tail:triangular} have modulus at least $3/4$,
so their inverse product has modulus at most $16/9<2$.  Also,
\[
 |C_0|<\frac12,\qquad |-1+h\nu|\le\frac98,
 \qquad |4-(6h+4)\nu|<5.
\]
Induction bounds $|C_h|$, through index $n-1$, by the coefficient of
$y^h$ in the nonnegative formal solution $Y(y)$ of
\begin{equation}\label{tail:majorant}
 Y=\frac12+3yY+10yY^2+2y^2Y^2+8y^2Y^3.
\end{equation}
At $y=1/64$, the right-hand side evaluated at $Y=1$ is
$1/2+13/64+10/4096<1$.  Iteration from $Y_0=0$ increases the
coefficients and stabilizes them successively to those of the unique
formal solution of \eqref{tail:majorant}.  Each iterate is at most
one at $1/64$.
Monotone convergence therefore gives
\[
 \sum_{h\ge0}[y^h]Y\,64^{-h}\le1,
 \qquad |C_h|\le[y^h]Y\le64^h.
\]
Equation~\eqref{tail:kappa-C} proves \eqref{tail:ordinary-bound}.
At $\nu=0$, this says $\kappa_1^\Lambda(b)=t$ and
$\kappa_r^\Lambda(b)=0$ for $r\ge2$.  At $t=0$, all $C_h$ vanish,
in agreement with \eqref{tail:t-zero}.
\end{proof}

The change of parameter gives the further bound
\begin{equation}\label{tail:u-ordinary-bound}
 |\kappa_r^\Lambda(b)|\le(r-1)!(128|u|)^{r-1}
 \quad\left(1\le r\le n,\ |u|\le\frac1{16n}\right),
\end{equation}
because $\nu=u/(1+u)$ then satisfies $|\nu|\le2|u|\le1/(8n)$.

\subsection{Norms of the transformed source}

For $G(z,b)=\sum_{\ell\ge0}G_\ell(b)z^\ell$ define
\[
 \|G\|_\rho=\sum_{\ell\ge0}\|G_\ell\|_1\rho^\ell
\]
when the sum is finite.  This norm is submultiplicative.  The basic
source estimate is
\begin{equation}\label{tail:six-bound}
 \|E_\ell\|_1,\ \|R_\ell\|_1\le6^\ell\qquad(\ell\ge1).
\end{equation}
To prove it, specialize to $x=y=1/2$.  The inverse relation for $T$
becomes $z=Te^{-T/2}$, giving $H=T/2$, $R=T/z$, and
$E=T'e^{-T}$.  Lagrange inversion yields
\begin{equation}\label{tail:diagonal}
 R_\ell(1/4)=\frac{(\ell+1)^\ell}
 {2^\ell(\ell+1)!},\qquad
 E_\ell(1/4)=\frac{(\ell-1)^\ell}{2^\ell\ell!}
 \quad(\ell\ge2).
\end{equation}
For the second identity, apply inversion to $E=-\partial_z(e^{-T})$.
The inequalities $n!\ge(n/e)^n$ and $e<3$ bound both expressions
by $(3/2)^\ell$ for $\ell\ge2$.  Nonnegativity and the degree
bound in \eqref{tail:source-coefficients} give
$P(1)\le4^\ell P(1/4)$ for each coefficient polynomial $P$.
This proves \eqref{tail:six-bound}, with the first-order cases
following from $R_1=(5+4b)/12$ and $E_1=0$.

\begin{proof}[Proof of Proposition~\ref{tail:quantitative-inputs}(ii)]
Lagrange inversion for \eqref{tail:change} gives
\begin{equation}\label{tail:alpha-Lagrange}
 \alpha_r=\frac1r[z^{r-1}]\frac{R_z}{R}R^{-r}
 =-\frac1r[z^r]R^{-r}.
\end{equation}
Taking $\rho=1/24$ in \eqref{tail:six-bound} gives
$\|R-1\|_\rho\le1/3$.  Expanding the inverse yields
$\|R^{-r}\|_\rho\le(3/2)^r$, so
\eqref{tail:alpha-Lagrange} gives $\|\alpha_r\|_1\le36^r/r$.

Put $J=1+zR_z/R$ and $G=\log E-\log J$.  The second Lagrange
formula is
\begin{equation}\label{tail:delta-Lagrange}
 \delta_r=\frac1r[z^{r-1}]G_zR^{-r}.
\end{equation}
Now take $\rho=1/48$.  Summing the geometric majorants in
\eqref{tail:six-bound} and their first two derivatives gives
\begin{gather*}
 \|R-1\|_\rho,\|E-1\|_\rho\le\frac17,
 \qquad \|R^{-1}\|_\rho,\|E^{-1}\|_\rho\le\frac76,\\
 \|R_z\|_\rho,\|E_z\|_\rho\le\frac{384}{49}<8,
 \qquad \|R_{zz}\|_\rho\le\frac{36864}{343}<108.
\end{gather*}
Thus $\|J-1\|_\rho<1/4$ and $\|J^{-1}\|_\rho\le4/3$.
Differentiating $G$ gives
\begin{align*}
 \|G_z\|_\rho
 &\le\|E_z/E\|_\rho+\|J^{-1}\|_\rho
 \left(\|R_z/R\|_\rho+
 \rho\|R_{zz}/R\|_\rho+
 \rho\|(R_z/R)^2\|_\rho\right)\\
 &\le\frac{28}{3}
 +\frac43\left(\frac{28}{3}+\frac{21}{8}+\frac{49}{27}\right)
 <48.
\end{align*}
Thus \eqref{tail:delta-Lagrange} implies
\[
 \|\delta_r\|_1\le
 \frac{48^{r-1}}r\,48\left(\frac76\right)^r
 =\frac{56^r}{r}.
\]
Products, inverses, and logarithms of series with constant term one
preserve the coefficient-degree bound $\deg_b[z^\ell]G\le\ell$.
In \eqref{tail:alpha-Lagrange} and \eqref{tail:delta-Lagrange},
differentiation increases the allowed degree by one.  It follows that
$\deg\alpha_r,\deg\delta_r\le r$.

For $|u|\le1/16$, \eqref{tail:conditional-cumulants} now gives
\[
 \|s_r\|_1\le(r-1)!56^r|u|^{r-1}
 \left((9/14)^r+|u|\right)
 \le(r-1)!56^r|u|^{r-1}.
\]
For $r\ge2$, apply $56^r\le56^{2r-2}=3136^{r-1}$.
For $r=1$, the identities $\alpha_1=c(b)$ and $\delta_1=-c(b)$
give
\[
 s_1=c(b)(1-u),\qquad
 \|s_1\|_1\le\frac34\left(1+\frac1{16}\right)
 =\frac{51}{64}<1.
\]
Together these bounds prove \eqref{tail:s-bound}.
\end{proof}

\subsection{Polynomial inputs and connected partitions}

\begin{proof}[Proof of Proposition~\ref{tail:quantitative-inputs}(iii)]
By multilinearity, it suffices to consider the inputs $b^{p_i}$.
Unitality makes the cumulant vanish when $a\ge2$ and an input is
constant; the case $a=1,p_1=0$ is immediate.  In the remaining
cases, put $N=\sum_i p_i\le n$ and partition $[N]$ into $a$
groups of sizes $p_i$, denoting this partition by $\pi_0$.
The product-cumulant formula gives
\begin{equation}\label{tail:product-cumulants}
 \kappa_\Lambda(b^{p_1},\ldots,b^{p_a})
 =\sum_{\substack{\tau\in\Pi([N])\\
                 \tau\vee\pi_0=\{[N]\}}}
   \prod_{B\in\tau}\kappa_{|B|}^\Lambda(b).
\end{equation}
To obtain this identity, expand the moments in
\eqref{tail:joint-definition} into partition cumulants.  For a fixed
$\tau$, the coefficient is the partition-lattice sum over partitions
of the $a$ groups coarser than the connections induced by $\tau$.
This sum equals one when all groups are connected and zero otherwise,
yielding the condition in \eqref{tail:product-cumulants}.

Let $\beta=64|\nu|$.  For $|u|\le1/(1024n)$,
\[
 |\nu|\le2|u|,\qquad
 \beta\le128|u|\le\frac1{8n}.
\]
Apply part (i) to each block in \eqref{tail:product-cumulants}.
Connectivity implies $N-|\tau|\ge a-1$.  For a fixed partition
$\tau$, the product $\prod_{B\in\tau}(|B|-1)!$ counts the permutations whose cycle
partition is $\tau$.  Enlarging the sum to all permutations with
this lower bound on cycle defect bounds \eqref{tail:product-cumulants}
in modulus by
\begin{equation}\label{tail:elementary-tail}
 \sum_{l\ge a-1}e_l(0,1,\ldots,N-1)\beta^l.
\end{equation}
Here $e_l$ denotes the elementary symmetric polynomial.  The identity
\[
 \sum_{\sigma\in S_N}\beta^{N-\#\operatorname{cycles}(\sigma)}
 =\prod_{i=0}^{N-1}(1+i\beta)
\]
for the symmetric group $S_N$ follows by inserting the last element
either as a singleton cycle or in one of the $N-1$ positions of an
existing cycle.

Put $h=a-1$ and $S=N(N-1)/2\le N^2/2$.  Since
$e_l(0,\ldots,N-1)\le S^l/l!$, we have
\begin{align*}
 \sum_{l\ge h}e_l(0,\ldots,N-1)\beta^l
 &\le\beta^h\frac{S^h}{h!}e^{\beta S}\\
 &\le h!\beta^h e^{2N+N/16}
 \le h!\beta^h8^n.
\end{align*}
For the middle step use $S^h/h!\le h!(N^h/h!)^2\le h!e^{2N}$
and $\beta S\le N^2/(16n)\le N/16$.
The last inequality uses $e^{33/16}<8$.  Indeed, the exponential
series gives $e<87/32$ and $e^{1/16}\le16/15$, hence
$e^{33/16}<7569/960<8$.  Using $\beta\le128|u|$ and summing the
monomial bounds with absolute coefficients proves
\eqref{tail:joint-bound}.
\end{proof}

\section{Polynomial inequalities at fixed Schur index}\label{supp:fixed-section}

We derive a family of kernels for Schur indices two through six
and prove the prefix inequalities that preserve log-concavity.
The parameter domains are partitioned explicitly, including the
exceptional prefixes for which the boundary and remainder terms
must be considered together.

For $d\in\mathbb Z_{\ge0}$, define
\[
 C_{k,d}(x,y)=[z^k]\prod_{i=0}^{d}
       \bigl(1+z(ix+(d-i)y)\bigr),\qquad
 A_{k,j}(d)=[u^j](1-u)C_{k,d}(1,u),
\]
and write $A_{k,j}(d)=\sum_NB_{k,j,N}\binom dN$.
Sequences are extended by zero outside their support.

\begin{theorem}\label{supp:fixed-indices}
For $2\le j\le6$ and every integer $k\ge2j$, the sequence
$(B_{k,j,N})_N$ is log-concave and has positive support
$k-1\le N\le2k$.
\end{theorem}

The cases $j=0,1$ are proved in \cite{BF}.  The proof here uses
the LC-positivity mechanism of Wang--Yeh \cite{WY}, applied to
explicit kernels derived from the original product.

\subsection{The source and its kernels}

Put $S_h(d)=e_h(0,1,\ldots,d)$.  Use
$i+(d-i)u=(1-u)i+du$ inside the product, then extract
$[u^j](1-u)$ to obtain
\begin{equation}\label{supp:source}
 A_{k,j}(d)=\sum_{s=0}^{j}(-1)^{j-s}\binom{k-s+1}{j-s}
       d^s\binom{d+1-k+s}{s}S_{k-s}(d).
\end{equation}
Indeed, after selecting the $k-s$ indices that contribute their
$(1-u)i$ terms, the remaining $s$ indices can be selected in
$\binom{d+1-k+s}{s}$ ways.  The other binomial coefficient in
\eqref{supp:source} comes from $(1-u)^{k-s+1}$.

For $f(d)=\sum_n b_n\binom dn$, let
$\mathcal Nf(t)=\sum_n b_nt^n$.  Multiplication and finite
differences give
\begin{equation}\label{supp:Newton-operators}
 \begin{aligned}
 \mathcal N(df)&=M\mathcal Nf,&M&=t(1+t)\partial_t+t,\\
 \mathcal NS_{h+1}&=L\mathcal NS_h,&L&=t(M+1).
 \end{aligned}
\end{equation}
The first identity follows from
$d\binom dn=n\binom dn+(n+1)\binom d{n+1}$.
For the second, use $\Delta S_{h+1}(d)=(d+1)S_h(d)$ and
$S_{h+1}(0)=0$.  These recurrences also prove
$\deg S_h\le2h$ and hence $\deg A_{k,j}\le2k$.

Set $m=k-j-1$ and $K=k-j=m+1$.  Define
\begin{equation}\label{supp:input-recurrence}
 R_0=1,\qquad
 R_{m+1}=\bigl(m+2+(m+3)t\bigr)R_m+t(1+t)R'_m.
\end{equation}
Starting from $\mathcal NS_1=t(1+t)$ gives
\begin{equation}\label{supp:input}
 \mathcal NS_K=t^K(1+t)R_m,
 \qquad R_m(t)=\sum_{a=0}^{m}w_at^a(1+t)^{m-a}.
\end{equation}
The polynomial $R_m$ has positive leading coefficient,
$R_m(0)=(m+1)!$, and $R_m(-1)=(-1)^m$.  Inductively, its
$m$ roots are simple and lie in $(-1,0)$: at consecutive roots
of $R_m$, the term $t(1+t)R'_m$ alternates in sign, and the
two endpoint signs give one new root in each interval cut out
by $-1$, the old roots, and $0$.  Degree counting proves that
these are all the roots of $R_{m+1}$.  Consequently
\[
 \sum_{a=0}^{m}w_az^a
       =(1-z)^mR_m\!\left(\frac{z}{1-z}\right)
\]
has only negative roots and positive leading coefficient.
Its coefficients are strictly positive and form a log-concave
sequence by Newton's inequalities.  They also satisfy
\begin{equation}\label{supp:weight-recurrence}
 w_a^{(m)}=(m+1+a)w_a^{(m-1)}+(m+1-a)w_{a-1}^{(m-1)},
 \qquad w_0^{(0)}=1.
\end{equation}

For independent $a,b$, put $k=a+b+j+1$ and define
\[
 \widetilde M=t(1+t)\partial_t+2a+b+1+(2a+2b+3)t,
 \qquad \widetilde L=t(\widetilde M+1).
\]
Our degree-$2j$ kernel is the following finite operator expression:
\begin{align}
 P^{(j)}_{a,b}(t)
 =j!\sum_{s=0}^{j}&\frac{(-1)^{j-s}}{s!}\binom{k-s+1}{j-s}
      \widetilde M^s
      \prod_{v=0}^{s-1}\bigl(\widetilde M-(k-s-1+v)\bigr)
      \widetilde L^{j-s}1\notag\\
 &=\sum_{h=0}^{2j}p_h^{(j)}(a,b)t^h.
       \label{supp:kernel-definition}
\end{align}
Operators act from right to left.  The factors involving only
$\widetilde M$ commute, but cannot be interchanged with
$\widetilde L$.  The coefficients $p_h^{(j)}$ are integer
polynomials: the binomial denominator in the source is cancelled
by $j!/s!$, and all remaining operations preserve integrality.
We suppress the superscript $j$ when it is fixed.

Conjugating \eqref{supp:Newton-operators} by
$t^{K+a}(1+t)^{b+1}$ and applying \eqref{supp:source} gives
\begin{equation}\label{supp:kernel-source}
 \begin{gathered}
 \sum_N B_{k,j,N}t^N=t^{k-1}\sum_{a=0}^{m}w_aF_a(t),\\
 F_a(t)=\frac1{j!}t^{a-j+1}(1+t)^{m-a+1}P^{(j)}_{a,m-a}(t).
 \end{gathered}
\end{equation}
Thus the output index $r$ below corresponds to $N=r+k-1$.
The coefficient matrix is
\begin{equation}\label{supp:kernel-entries}
 K_{r,a}=[t^r]F_a(t)
   =\frac1{j!}\sum_{h=0}^{2j}p_h(a,m-a)
             \binom{m-a+1}{r-a+j-1-h}.
\end{equation}
Binomial coefficients are zero outside their natural range when
their upper index is a nonnegative integer.

\subsection{Columns and the exact range of the kernel argument}

The thresholds used in the proof are recorded here:
\begin{equation}\label{supp:thresholds}
\begin{array}{c|rrrrr}
 j&2&3&4&5&6\\\hline
 M_j\text{ (minimum }m\text{)}&1&2&3&5&7\\
 q_j\text{ (minimum interior }q\text{)}&1&1&1&1&2\\
 t_j\text{ (number of upper centers)}&4&6&8&10&13\\
 \text{remaining degrees }k&\text{none}&\text{none}&\text{none}&10&12,13
\end{array}
\end{equation}
The column argument applies to $a,b\ge0$, $a+b\ge M_j$.
This domain has the disjoint decomposition
\begin{equation}\label{supp:column-domains}
 \begin{aligned}
 &(a,b)=(A+M_j,B),\qquad A,B\ge0,\\
 &(a,b)=(i,B+M_j-i),\qquad 0\le i<M_j,\quad B\ge0.
 \end{aligned}
\end{equation}
For $j=3,4,5,6$, expansion of \eqref{supp:kernel-definition}
after each substitution in \eqref{supp:column-domains} gives
\[
 p_h\in\mathbb Z_{\ge0}[A,B]\quad(0\le h\le2j),\qquad
 p_h^2-p_{h-1}p_{h+1}\in\mathbb Z_{\ge0}[A,B]
 \quad(1\le h<2j),
\]
with the unused variable omitted in a one-variable domain.
For $j=2$, these coefficient inequalities hold directly in
$\mathbb Z_{\ge0}[a,b]$.

The first nonzero coefficient of $P^{(j)}_{a,b}$ is normally
$h=\max(0,j-1-a)$, and every coefficient from that index through
$2j$ is positive.  There is one exception: for $j=4$, $a=0$,
$b=m=3$, the coefficient with $h=3$ also vanishes.  Explicitly,
\[
 p_3^{(4)}(0,b)=4(b-3)(b+2)(b+4)^3(b^2+5b+5),
\]
and $p_4^{(4)}(0,b),\ldots,p_8^{(4)}(0,b)$ are positive for
$b\ge3$.  In particular every apparent negative power in
\eqref{supp:kernel-source} cancels.

Multiplication by $1+t$ preserves nonnegative log-concavity and
interval support.  For $d_i=c_i+c_{i-1}$, its Tur\'an difference is
\[
 (c_i^2-c_{i-1}c_{i+1})+(c_{i-1}^2-c_{i-2}c_i)
       +(c_ic_{i-1}-c_{i-2}c_{i+1})\ge0;
\]
the last term is nonnegative by decreasing consecutive ratios,
including zero extension.  It follows that every column $F_a$
is log-concave and has positive support
\begin{equation}\label{supp:column-support}
 \max(0,a-j+1)\le r\le m+j+2,
\end{equation}
except that the column $(j,m,a)=(4,3,0)$ starts at $r=1$.
The right endpoint is common to all columns.  The column $F_0$
has full support starting at zero except in that one case, where
$F_1$ has full support starting at zero.

\subsection{The centered prefix criterion}

For $0\le S\le2m$ and
$\max(0,S-m)\le a\le\lfloor S/2\rfloor$, define
\begin{equation}\label{supp:prefix}
 \mathcal P_{r,S,a}=
 \sum_{c=a}^{S-a}
   \bigl(K_{r,c}K_{r,S-c}-K_{r-1,c}K_{r+1,S-c}\bigr).
\end{equation}
The sum contains each off-diagonal pair in both orders and a
middle diagonal once.

\begin{lemma}\label{supp:Abel}
If every prefix in \eqref{supp:prefix} is nonnegative at every
interior output index, the kernel maps every nonnegative
log-concave input with interval support to a log-concave sequence
with interval support.
\end{lemma}

\begin{proof}
Fix $S$ and combine the off-diagonal pair with indices $a,S-a$
into
\[
 T_a=2K_{r,a}K_{r,S-a}-K_{r-1,a}K_{r+1,S-a}
                            -K_{r-1,S-a}K_{r+1,a}.
\]
If $a=S/2$, retain its single diagonal instead.  For
$l=\max(0,S-m)$ and $h=\lfloor S/2\rfloor$, put
$C_a=\sum_{c=a}^{h}T_c=\mathcal P_{r,S,a}$ and
$z_a=w_aw_{S-a}$.  Log-concavity of the input makes $z_a$
nondecreasing toward the center, also when initial entries are
zero.  Abel summation gives
\[
 \sum_{a=l}^{h}T_az_a
       =C_lz_l+\sum_{a=l+1}^{h}C_a(z_a-z_{a-1})\ge0.
\]
Summing over $S$ proves each output Tur\'an inequality.  The union
of the nonempty column supports is an interval because they have
a common right endpoint.  This is the LC-positivity argument of
\cite{WY} in the present notation.
\end{proof}

Diagonal prefixes have already been proved by column
log-concavity.  We next prove every off-diagonal prefix, first in
the interior and then at the upper centers.

\subsection{A common-binomial and Gosper identity}

For the interior, set
\begin{equation}\label{supp:interior-coordinates}
 \rho=r+j-2\ge j-1,\qquad q=m-\rho\ge q_j,
 \qquad D_q=(q+1)^{\overline{2j}}.
\end{equation}
Here $v^{\overline h}=v(v+1)\cdots(v+h-1)$ is a rising
factorial.  The interior range is
$1\le r\le m-j+2-q_j$.  In all adjacent-row expressions,
$m=\rho+q$ is fixed.  Define
\begin{equation}\label{supp:H}
 H_v(a)=\sum_{h=0}^{2j}p_h(a,m-a)
       \fall{v-a+1}{h}(m-v+1+h)^{\overline{2j-h}}.
\end{equation}
Complementing the binomial coefficient in
\eqref{supp:kernel-entries} gives
\begin{equation}\label{supp:common-binomial}
 K_{r,a}=\frac1{j!}\sum_{h=0}^{2j}p_h(a,m-a)
                    \binom{m-a+1}{q+h}
       =\frac{\binom{m-a+1}{q}H_\rho(a)}{j!D_q}.
\end{equation}
If the common binomial is zero, all binomials in the sum are
zero; otherwise the falling factors in $H_\rho$ account for
terminal zeros.  This is an identity without division by a
possibly zero binomial coefficient.

Put $b=S-a$ and $M_a=m-a+1$, and define
\begin{align}
 W(a)&=\frac{\binom{M_a}{q}\binom{M_b}{q-1}}{q(q+1)},
 &G(a)&=\binom{M_a+1}{q+1}\binom{M_b}{q},\notag\\
 \mathsf A(a)&=(\rho+q-a+2)(\rho-b+2),
 &\mathsf B(a)&=(\rho-a+1)(\rho+q-b+2).
       \label{supp:weights}
\end{align}
Set
\begin{equation}\label{supp:Phi}
 \begin{aligned}
 \Phi(a)=(q+1)\bigg[&(\rho-b+2)H_\rho(a)H_\rho(b)\\
   &-\frac{q+2j}{q+2j+1}(\rho-a+1)
               H_{\rho-1}(a)H_{\rho+1}(b)\bigg].
 \end{aligned}
\end{equation}
Direct binomial multiplication identities give
\begin{align}
 (j!)^2D_q^2\bigl(K_{r,a}K_{r,b}-K_{r-1,a}K_{r+1,b}\bigr)
     &=W(a)\Phi(a),\label{supp:scaled-summand}\\
 W(a)\mathsf A(a)=G(a),\qquad W(a)\mathsf B(a)&=G(a+1).
       \label{supp:weight-telescope}
\end{align}
The normalization uses
$D_q^2/(D_{q+1}D_{q-1})=(q+1)(q+2j)/(q(q+2j+1))$.
All $W(a)$ are nonnegative on the legal domain, including $q=1$
and binomial endpoints.

For $j=2,3,4,5,6$, set
\begin{equation}\label{supp:Gosper-denominator}
 E_j=2^{2j}(q+2j+1)\prod_{v=0}^{2j}(2q+2v+1).
\end{equation}
Define $U$ and the $a$-independent remainder $R$ by the triangular
elimination below, and put $U_{\rm num}=E_jU$ and
$R_{\rm num}=E_jR$.  These are integer polynomials.  There are
integer polynomials $V,R_{\rm core}$ satisfying
\begin{equation}\label{supp:Gosper}
 \Phi(a)=\mathsf A(a)U(a)-\mathsf B(a)U(a+1)+R,
\end{equation}
\begin{equation}\label{supp:reflection}
 \begin{aligned}
 (\rho-a+2)U_{\rm num}(a)&-(\rho-b+1)U_{\rm num}(b+1)\\
        &=(b-a+1)V(a,q,\rho,S),
 \end{aligned}
\end{equation}
\begin{equation}\label{supp:remainder-factor}
 R_{\rm num}=\Lambda_j R_{\rm core},
\end{equation}
where the suppressed arguments $q,\rho,S$ are unchanged and
\begin{equation}\label{supp:positive-factor}
 \Lambda_j=q(q+1)^2\prod_{h=2}^{2j+1}(q+h)
       (2\rho-S+3)(2q+2\rho-S+4)(2q+2\rho-S+5).
\end{equation}
The degree of $U_{\rm num}$ in $a$ is at most $4j$; both
remainder polynomials are independent of $a$.

To specify $U$ and $R$, expand the finite source, for which
$\deg_a\Phi\le4j+1$.
Starting with $h=4j$ and proceeding down to zero, let $c$ be the
current coefficient of $a^{h+1}$ and subtract
\[
 \frac{c}{2q+h+1}
       \bigl(\mathsf A(a)a^h-\mathsf B(a)(a+1)^h\bigr).
\]
The leading coefficient of the parenthesized polynomial is
$2q+h+1$, so this removes its highest remaining degree.  The
collected coefficients define $U$, and the constant left in $a$
defines $R$.  Multiplication by $E_j$ gives the integer witnesses
above.  Polynomial division gives $V$ and $R_{\rm core}$ through
\eqref{supp:reflection} and \eqref{supp:remainder-factor}.
This specifies every term of the identities from
\eqref{supp:kernel-definition}.

Summing \eqref{supp:scaled-summand} and telescoping gives
\begin{equation}\label{supp:telescoped}
 (j!)^2D_q^2\mathcal P_{r,S,a}
 =G(a)U(a)-G(b+1)U(b+1)+R\sum_{c=a}^{b}W(c).
\end{equation}
If $b\le\rho+1$, its boundary part is
\begin{equation}\label{supp:boundary}
 \frac{G(a)(b-a+1)}{E_j(\rho-a+2)}V(a,q,\rho,S).
\end{equation}
Indeed $G(b+1)(\rho-a+2)=G(a)(\rho-b+1)$, also when
$b=\rho+1$ and $G(b+1)=0$.  If $b=\rho+2$, both boundary
weights vanish.  The denominators are positive, and
$\Lambda_j>0$ whenever $q\ge1$ and $S\le2\rho+2$.
It remains to establish the indicated signs of $R_{\rm core}$
and $V$, or to retain an entire exceptional prefix.

\subsection{Exhaustive interior domains}\label{supp:domains}

The bound \eqref{supp:column-support} implies
$K_{r,a}=0$ for $a>r+j-1=\rho+1$.  A central product therefore
requires $a,b\le\rho+1$, and an adjacent-row product requires
$a\le\rho$, $b\le\rho+2$.  Entire traces vanish when
$S>2\rho+2$.  If an outer endpoint $b>\rho+2$, all three
entries at that endpoint vanish, so both ordered outer terms
can be deleted.  Repeating leaves an empty prefix, a diagonal,
or a non-diagonal interval satisfying
\begin{equation}\label{supp:active-domain}
 \begin{gathered}
 0\le a<b\le m,\qquad \rho\ge j-1,\qquad q\ge q_j,\\
 S=a+b\le2\rho+2,\qquad b\le\rho+2.
 \end{gathered}
\end{equation}

For $j=2$, the following four domains cover every such interval.
All free variables are nonnegative integers, and the last column
lists the polynomial signs required in \eqref{supp:telescoped}.
\begingroup\setlength{\arraycolsep}{3pt}
\begin{equation}\label{supp:j2-domains}
\begin{array}{c|cccc|c}
 \text{domain}&a&q&\rho&S&\text{signs}\\\hline
 b\le\rho&L&Q+1&L+D+Z+1&2L+D+1&R_{\rm core},V\\
 b=\rho+1,\ b-a\ge2&L&Q+1&L+D+1&2L+D+2&R_{\rm core},V\\
 b=\rho+1,\ b-a=1&L+1&Q+1&L+1&2L+3&R_{\rm core},V\\
 b=\rho+2&L&Q+2&L+D&2L+D+2&R_{\rm core}
\end{array}
\end{equation}
\endgroup
The last row may include additional points with $\rho=0$; its
coefficient certificate is valid there as well.  The first row
parametrizes the positive gap and the excess $\rho-b$.  In the
second and third rows the gap is at least two or is one; in the
last row the trace bound forces a gap at least two and $b\le m$
forces $q\ge2$.  The required polynomials after these substitutions have
nonnegative integer coefficients.

For $3\le j\le6$, put $\ell=j-1$ and define the exceptional
intervals
\begin{equation}\label{supp:exceptions}
 \mathcal E_3=\varnothing,\quad
 \mathcal E_4=\{(0,1)\},\quad
 \mathcal E_5=\mathcal E_6=\{(0,1),(0,2)\}.
\end{equation}
Table~\ref{tab:interior-domains} specifies all other domains.  In the
first five rows take $q=Q+q_j$; in the last two take $q=Q+2$.
The integers $L,D,Z,Q$ are nonnegative.

\begin{table}[htbp]
\caption{Interior parameter domains for $3\le j\le6$.}
\label{tab:interior-domains}
\centering\footnotesize\setlength{\tabcolsep}{2pt}
\begin{tabular}{@{}ccccc@{}}
\toprule
 Domain & $a$ & $\rho$ & $S$ & Signs\\
\midrule
 $b\le\rho$, $a\ge\ell-1$
  &$L+\ell-1$&$L+D+Z+\ell$&$2L+D+2\ell-1$&$R_{\rm core},V$\\
 $b\le\rho$, $0\le i\le\ell-2$, $b\ge\ell$
  &$i$&$D+Z+\ell$&$D+\ell+i$&$R_{\rm core},V$\\
 $0\le i\le\ell-2$, $1\le h\le\ell-1-i$
  &$i$&$Z+\ell$&$2i+h$&$R_{\rm core},V$\\
 $b=\rho+1$, $1\le h\le\ell$
  &$L+\ell+1-h$&$L+\ell$&$2L+2\ell+2-h$&$R_{\rm core},V$\\
 $b=\rho+1$, $b-a\ge\ell+1$
  &$L$&$L+D+\ell$&$2L+D+\ell+1$&$R_{\rm core},V$\\
 $b=\rho+2$, $2\le h\le\ell+1$
  &$L+\ell+2-h$&$L+\ell$&$2L+2\ell+4-h$&$R_{\rm core}$\\
 $b=\rho+2$, $b-a\ge\ell+2$
  &$L$&$L+D+\ell$&$2L+D+\ell+2$&$R_{\rm core}$\\
\bottomrule
\end{tabular}
\end{table}

In the third row omit $(i,i+h)\in\mathcal E_j$.  In the
fourth and sixth rows, $h=b-a$ is the gap.  After each indicated substitution, the required polynomials have
nonnegative monomial expansions in the displayed free variables.
Since those variables range over nonnegative integers, the
remainder and boundary terms in \eqref{supp:telescoped} are
nonnegative throughout each domain.

To check exhaustion, first suppose $b\le\rho$.  Either
$a\ge\ell-1$, or $a\le\ell-2$ and $b\ge\ell$, or
$0\le a<b\le\ell-1$.  These are the first three rows, with
the listed exceptions removed.  If $b=\rho+1$, then
$a+h\ge\ell+1$: either $1\le h\le\ell$, with
$a\ge\ell+1-h$, or $h\ge\ell+1$.  If $b=\rho+2$, the
trace bound forces $h\ge2$, while
$a+h\ge\ell+2$ and $q\ge2$.  Distinguishing
$2\le h\le\ell+1$ from $h\ge\ell+2$ gives the last two
rows.  These subdivisions are disjoint and exhaust
\eqref{supp:active-domain} except for \eqref{supp:exceptions}.

\subsection{Complete exceptional prefixes}\label{supp:corners}

The separated remainder does not have the required sign at the
exceptional intervals, so we use their complete ordered sums.
Let $\Phi_c^{(S)}$ denote \eqref{supp:Phi} at $a=c$ and trace
$S$, and define
\begin{align}
 Y_1={}&(\rho+2)(q+2j+1)\Phi_0^{(1)}
             +(\rho+1)(q+2j+1)\Phi_1^{(1)},\label{supp:Y1}\\
 Y_2={}&(\rho+2)(m+1)(q+2j+1)\Phi_0^{(2)}\notag\\
       &+(\rho+2)m(q+2j+1)\Phi_1^{(2)}
             +\rho(m+1)(q+2j+1)\Phi_2^{(2)}.
             \label{supp:Y2}
\end{align}
The denominators in $\Phi$ cancel, so these are polynomials.
At trace one, $W(1)/W(0)=(\rho+1)/(\rho+2)$; at trace two,
the corresponding ratios are $m/(m+1)$ and $\rho/(\rho+2)$.
Therefore
\begin{align}
 \mathcal P_{r,1,0}
  &=\frac{W(0)Y_1}{(j!)^2D_q^2(\rho+2)(q+2j+1)},
       \label{supp:corner-one}\\
 \mathcal P_{r,2,0}
  &=\frac{W(0)Y_2}{(j!)^2D_q^2(\rho+2)(m+1)(q+2j+1)}.
       \label{supp:corner-two}
\end{align}
Each $W(0)$ is evaluated at its own trace and is positive.
The trace-two prefix includes its middle diagonal exactly once.
The required polynomial signs are
\begin{equation}\label{supp:corner-inventory}
\begin{array}{c|c|c}
 j&(q,\rho)&\text{polynomials in }\mathbb Q_{\ge0}[Q,Z]\\\hline
 4&(Q+1,Z+3)&Y_1\\
 5&(Q+1,Z+4)&Y_1,Y_2\\
 6&(Q+2,Z+5)&Y_1,Y_2.
\end{array}
\end{equation}
They follow by expanding \eqref{supp:Y1}--\eqref{supp:Y2} after
the indicated substitutions.  Equations
\eqref{supp:corner-one}--\eqref{supp:corner-two} then prove each
exceptional prefix throughout its domain.  Together with the
preceding domain tables and support trimming, this proves all
interior prefixes.

\subsection{All upper centers}\label{supp:upper-centers}

For $0\le\delta\le t_j+1$, let
\begin{equation}\label{supp:upper-rows}
 T_\delta(a,m)=K_{m+j+2-\delta,a}
   =\frac1{j!}\sum_{h=0}^{2j}p_h(a,m-a)
             \binom{m-a+1}{\delta+h-2j}.
\end{equation}
A negative lower binomial index contributes zero; for a fixed
nonnegative lower index the binomial is its falling-factorial
polynomial.  This agrees with the zero-extended convention on
the nonnegative integer upper indices used here.

Write $m=L+D+Y\ge M_j$, where $L,D,Y\ge0$.  At each
$1\le\delta\le t_j$, the complete ordered prefix is
\begin{equation}\label{supp:upper-prefix}
 \begin{aligned}
 \mathcal W_\delta(L,D,Y)=\sum_{i=0}^{D}\bigl[&
 T_\delta(L+i,m)T_\delta(L+D-i,m)\\
 &-T_{\delta+1}(L+i,m)T_{\delta-1}(L+D-i,m)\bigr].
 \end{aligned}
\end{equation}
For $j=2$ and $1\le\delta\le4$, this polynomial has the expansion
\[
 \mathcal W_\delta(L,D,Y)
      =\sum_{l,h,y}c_{l,h,y}^{(\delta)}L^l\binom DhY^y,
       \qquad c_{l,h,y}^{(\delta)}\ge0,
\]
on the entire nonnegative domain.

For $j=3,4,5,6$, diagonal prefixes $D=0$ have already been
handled.  All $D\ge1$ prefixes lie in exactly one of the
following domains, where $d$ ranges over $1\le d<M_j$ and
$l$ over $0\le l<M_j-d$.  Every free parameter is a
nonnegative integer.
\begin{equation}\label{supp:upper-domains}
\begin{array}{c|c}
 \text{substitution}&\text{free parameters}\\\hline
 D=Q+M_j&L,Q,Y\\
 D=d,\quad L=X+M_j-d&X,Y\\
 D=d,\quad L=l,\quad Y=Z+M_j-d-l&Z.
\end{array}
\end{equation}
First distinguish $D\ge M_j$ from $D=d<M_j$, and then
$L\ge M_j-d$ from $L<M_j-d$.  In the last case $m\ge M_j$
forces exactly the stated lower bound on $Y$.  This proves
disjointness and exhaustion.

In the first domain, expand in the basis
$L^l\binom QhY^y$; in the remaining domains, use ordinary
monomials in the displayed free variables.  The coefficients are
nonnegative rational numbers in each case.  These expansions give
$\mathcal W_\delta\ge0$, since every basis element is nonnegative
at the stated integer parameters.

For completeness, the polynomials in \eqref{supp:upper-prefix}
are determined by the power-sum identities
\begin{equation}\label{supp:power-sums}
 \sum_{i=0}^{D}i^p
  =\sum_{h=0}^{p}\genfrac{\{}{\}}{0pt}{}{p}{h}h!
          \binom{D+1}{h+1}
  =\frac{B_{p+1}(D+1)-B_{p+1}(0)}{p+1},
\end{equation}
where the braces denote Stirling numbers of the second kind and
$B_n$ are ordinary Bernoulli polynomials.  Thus both the power expansion and the mixed binomial expansion
are determined by the finite kernel formula
\eqref{supp:kernel-definition}.

\subsection{Completion and the finite degrees}\label{supp:completion}

For $m\ge M_j$, the interior proof covers
$1\le r\le m-j+2-q_j$.  Since
$t_j=2j-1+q_j$, the upper proof covers the immediately following
range
\[
 m-j+3-q_j\le r\le m+j+1.
\]
These are every nontrivial output center; the endpoints
$r=0,m+j+2$ are automatic under zero extension.  Diagonal
prefixes, support trimming, and all exceptional intervals have
been included.  Lemma~\ref{supp:Abel} therefore applies to the
actual strictly positive log-concave weights in
\eqref{supp:input}.  By \eqref{supp:column-support}, at least
one column is positive on all of $[0,m+j+2]$.  The output has
that exact support.  Multiplication by $t^{k-1}$ in
\eqref{supp:kernel-source} gives the positive support
$[k-1,2k]$ and proves log-concavity for every degree covered by
the kernel argument.

The only pairs below the kernel threshold are $(j,k)=(5,10)$,
$(6,12)$, and $(6,13)$.  In each case $j=\lfloor k/2\rfloor$.
Theorem~\ref{near:balanced}, followed by
Theorem~\ref{bridge:full}, proves their Newton log-concavity and
indeed the stronger factor \eqref{intro:factor}.  This completes
the proof of Theorem~\ref{supp:fixed-indices}.

\end{appendices}
\backmatter
\section*{Statements and Declarations}

\subsection*{Data availability}
The normalized coefficient rows for the finite central ranges in
Proposition~\ref{comp:central} are provided in Online Resource~1.
The polynomial witnesses and coefficient expansions for the fixed-index
inequalities in Appendix~\ref{supp:fixed-section} are provided in
Online Resource~2.
The generating identities, recurrences and parameter domains are
specified in the article.

\subsection*{Use of generative artificial intelligence}

This study was conducted by researchers in the field of AI. The authors are not mathematics specialists. The study was completed entirely using OpenAI's ChatGPT, including the organization of research ideas, literature search and synthesis, assistance with proof derivation, organization of computational assertions, preparation of verification code, and drafting and editing of the manuscript. To reduce factual, computational, and derivational errors in the generated content, the authors used Tencent WorkBuddy to invoke the DeepSeek model to verify key content, including the main computational assertions, mathematical reduction steps, executable checks, and some literature leads.

 Content generated by AI systems may still contain inaccuracies or incompleteness; responsibility for the content rests with the human authors. This study was conducted purely out of personal interest, does not represent the position of the authors' institutions, and does not constitute professional mathematical advice. Readers should regard this study as an exploratory work completed with AI assistance and should independently verify its key conclusions.


\end{document}